\input amstex
\magnification 1200
\vcorrection{-12mm}
\documentstyle{amsppt}
\NoBlackBoxes

\def\eps{\varepsilon}
\def\Irr{\operatorname{Irr}}
\def\UC{\Omega}                   

\def\codim{\operatorname{codim}}
\def\rk{\operatorname{rk}}
\def\cn{\operatorname{cn}}
\def\ecn{\operatorname{ecn}}
\def\GF#1{\Bbb F_{#1}}


\def\refAH         {1}
\def\refAW         {2}
\def\refBelkale    {3}
\def\refDickson    {4}
\def\refEnnolaConj {5}
\def\refEnnola     {6}
\def\refGAP        {7}
\def\refGeck       {8}
\def\refGordeev    {9}
\def\refKarni     {10}
\def\refLev       {11}
\def\refLSh       {12}
\def\refMM        {13}
\def\refOrJKTR    {14}
\def\refSF        {15}
\def\refWallG     {16}
\def\refWallSL    {17}


\def\eqSCF        {1}
\def\eqDet        {2}
\def\eqRk         {3}
\def\eqDeltaL     {4}
\def\eqEV         {5}
\def\eqSumOverIrr {6}

\def\eqCSF{\eqSCF}


\def\propEET     {1.1}
\def\propSUinv   {1.2}
\def\thGU        {1.3}
\def\remGU       {1.4}
\def\remGUbis    {1.5}
\def\thSU        {1.6}
\def\corPGU      {1.7}
\def\corPSU      {1.8}
\def\corECN      {1.9}
\def\remKarniLev {1.10}

\def\lemQone    {2.1}  
\def\lemQtwo    {2.2}

\def\propSmallqOne   {4.1}
\def\propSmallqTwo   {4.2}
\def\propSmallqThree {4.3}
\def\propSmallqSU    {4.4}
\def \corSmallqSU    {4.5}

\def\thGUtwo      {5.1}
\def\propSUtwoInv {5.2}
\def\remSUtwo     {5.3}
\def\thSUtwo      {5.4}
\def\remSUtwoF    {5.5}
\def\corPSUtwo    {5.6}


\def\sectDetRk      {1.2}
\def\sectConvent    {1.3}
\def\sectCCG        {1.4}
\def\sectCCS        {1.5}
\def\sectNotation   {1.6}
\def\sectMainTh     {1.7}
\def\sectCN         {1.8}

\def\sectG           {2}
\def\sectCT          {2.1}
\def\sectSCF         {2.2}
\def\sectTP          {2.3}
\def\sectQthreeQfour {2.4}
\def\sectTPtwo       {2.5}
\def\sectQP          {2.6}

\def\sectS          {3}
\def\sectCTS        {3.1}
\def\sectSCS        {3.2}
\def\sectProofS     {3.3}

\def\sectSmallq      {4}
\def\sectSmallqG     {4.1}
\def\sectSmallqS     {4.2}

\def\sectNTwo       {5}
\def\sectNTwoG      {5.1}
\def\sectNTwoCCS    {5.2}
\def\sectNTwoS      {5.3}
\def\sectNTwoP      {5.4}

\def\sectNtwo{\sectNTwo}
\def\sectNtwoP{\sectNTwoP}

\def\tabCC         {1}
\def\tabExcept     {2}
\def\tabCT         {3}

\def\tabTripleProd {4}
\def\tabTPExcept   {5}

\def\tabSLthreeQtwo   {6.1}
\def\tabSUthreeQthree {6.2}
\def\tabSUthreeQfour  {6.3}

\def\tabSU         {7}
\def\tabGUtwo      {8}
\def\tabSUtwoSLtwo {9}
\def\tabSUtwo     {10}
\def\tabPSUtwo    {11}


\rightheadtext{Products of conjugacy classes in $GU(3,q^2)$ and $SU(3,q^2)$}


\topmatter
\title     Products of conjugacy classes in finite unitary
           groups $GU(3,q^2)$ and $SU(3,q^2)$
\endtitle
\author    S.Yu.~Orevkov
\endauthor
\address   Institut des Math\'ematiques de Toulouse, UPS,
           118 route de Narbonne, 31062 Toulouse, France
\endaddress
\endtopmatter

\document

\head 1. Introduction and statement of main results
\endhead

\subhead 1.1. Introduction
\endsubhead
We study here the following problem (the {\it Class Product Problem}).
Let $c_1,\dots,c_m$ be conjugacy classes in a given group.
Does the unity of the group belong to their product?
For the usual unitary group $SU(n)$, this problem is completely solved
in [\refAW] and [\refBelkale].
Various partial cases of the class product problem (in particular, estimates for the
covering number) for many groups were studied by many authors 
see, e.~g., [\refAH, \refGordeev, \refLev, \refLSh] and numerous references therein.

In this paper we give
a complete solution to the class product problem for the finite unitary groups $GU(3,q^2)$
and $SU(3,q^2)$, see \S\sectMainTh\ for precise statements.
Due to Ennola duality (see \S\sectConvent), as a by-product, we obtain a solution
for the groups $GL(3,q)$, $SL(3,q)$.
For the sake of completeness, we also give in \S\sectNtwo\ a solution for the groups
$GL(2,q)$, $GU(2,q^2)$ and $SU(2,q^2)\!\cong\! SL(2,q)$.
A solution for corresponding projective groups $PGL$, $PSL$, $PGU$ and $PSU$
easily follows.

Our interest to the class product problem in all kinds of unitary groups is motivated
by the study of braid monodromy of plane algebraic curves (see [\refOrJKTR]).

As in [\refAH, \refLSh], the main tool used here for solving the class product problem
is Burnside's formula for the {\it structure constants} via the character table.
Namely, for a finite group $\Gamma$ and its elements $x_1,\dots,x_m$, we denote the number
of $m$-tuples $(y_1,\dots,y_m)$ such that $y_i$ is a conjugate of $x_i$ in $\Gamma$ and
$y_1\dots y_m=e$ by $N_\Gamma(x_1,\dots,x_m)$. Then Burnside's formula
(see, e.~g., [\refMM; Th.~I-5.8] or [\refAH; Ch.~1, 10.1]) reads as
$$
   N_\Gamma(x_1,\dots,x_m) = {|x_1^\Gamma|\cdot\dots\cdot|x_m^\Gamma|\over|\Gamma|}
    \sum_{\chi\in\Irr(\Gamma)} {\chi(x_1)\dots\chi(x_m)\over \chi(1)^{m-2}}
                                                                         \eqno(\eqSCF)
$$
where $\Irr(\Gamma)$ is the set of irreducible characters of $\Gamma$ and $x^\Gamma$ denotes the
conjugacy class of $x$ in $\Gamma$. We denote the sum in the right hand side of (\eqSCF)
by $\bar N_\Gamma(x_1,\dots,x_m)$.

We use the character tables from [\refEnnola] ($GU/GL$) and [\refSF, \refGAP] ($SU/SL$).

\subhead Acknowledgment
\endsubhead
I am grateful to M.~Geck, A.~A.~Klyachko, N.~A.~Vavilov and I.~A.~Vedenova for useful
advises and discussions.


\subhead\sectDetRk. Determinant Relation and Rank Condition
\endsubhead
If $\Gamma$ is a subgroup of $GL(n,K)$ over any commutative field $K$ and
$A_1,\dots,A_m\in \Gamma$ are such that $I\in A_1^\Gamma\dots A_m^\Gamma$, then
an evident restriction is the {\it determinant relation}
$$
  \det(A_1)\cdot\dots\cdot\det(A_m)=1		\eqno(\eqDet)
$$
Another evident restriction which takes place for any field, is the {\it rank condition}: if
$\lambda_1\dots\lambda_m=1$, then
$$
   \rk(A_j-\lambda_j I) \le \sum_{i\ne j}\rk(A_i - \lambda_i I)
   \qquad\text{for any $j=1,\dots,m$}
                                                                   \eqno(\eqRk)
$$
($I$ is the identity matrix).
Indeed, if we denote the $\lambda_i$-eigenspace of $A_i$ by $V_i$, then
$\bigcap_{i\ne j}V_i\subset V_j$, thus
$\codim V_j\le\codim\bigcap_{i\ne j} V_i\le\sum_{i\ne j}\codim V_i$.
When $m>n$, this condition is always satisfied for
any $m$-tuple of non-scalar matrices.

One more general restriction (see Case  $(viii)$ in Theorem \thGU(a)) is

\proclaim{ Proposition \propEET } Let $K$ be a perfect field and $A\sim B\in GL(3,K)$.
If $A$ does not have eigenvalues in $K$, then
$A^{-1}B\ne\left(\smallmatrix1&0&0\\1&1&0\\0&0&1\endsmallmatrix\right)$.
\endproclaim

\demo{ Proof }
Suppose the contrary. Let $V$ be the eigenspace of $A^{-1}B$.
Then $A|_V=B|_V$. Since $A$ has no eigenvalues in $K$, we have $A(V)\ne V$.
Let $e_2\in V\cap A(V)$, $e_1=A^{-1}(e_2)$, and $e_3=A(e_2)$.
Then $B(e_1)=A(e_1)=e_2$ and $B(e_2)=A(e_2)=e_3$. Thus, 
$A$ and $B$ take the canonical form
in the same basis $(e_1,e_2,e_3)$. Since  $A\sim B$, this implies $A=B$. Contradiction.
\qed
\enddemo

It happens (see Theorem \thGU\ in \S\sectMainTh) that in the case of $GL(3,q)$, $q\ne2$,
there are no other restrictions on
$A_1,\dots,A_m$. In the case of $GU(3,q^2)$, there are much more restrictions (see the lines
in Table \tabExcept\ not marked by the asterisk).
An interesting question is to generalize them for any field and for any dimension.


\subhead \sectConvent. Ennola duality and the sign convention
\endsubhead
Throughout the paper, $q$ is a prime power and
$GU$ (resp. $SU$, $GL$, $SL$) is an abbreviation of
$GU(3,q^2)$ (resp. $SU(3,q^2)$, $GL(3,q)$, $SL(3,q)$)
except \S\sectNTwo\ where the same convention is used with $3$ replaced by $2$.

Ennola [\refEnnola] observed that the character tables of groups $GU(n,q^2)$ and
$GL(n,q)$ are obtained from each other by changing the sign of $q$.
The same is true for $SU(n,q^2)$ and $SL(n,q)$.
Since the character table is our main tool, it is not surprising that all
computations are almost the same for $GU/SU$ and $GL/SL$.
So, throughout the paper (except \S\sectSmallq\ and \S\sectNTwoS),
we use the following {\it sign convention}: if a symbol $\pm$ or $\mp$ occurs in
a formula, then the upper sign corresponds to the case of $GU$ (resp. $SU$, $PSU$)
and the lower sign corresponds to the case of $GL$ (resp. $SL$, $PSL$).
Throughout the paper (except \S\sectSmallq\ and \S\sectNTwoS),
$G$ (resp. $S$; $PG$; $PS$) stands for
$GU$ or $GL$ (resp $SU$ or $SL$; $PGU$ or $PGL$; $PSU$ or $PSL$) and we set
$$
    \delta_L = {1\mp1\over2} = \cases 1, &G=GL,\\ 0, &G=GU.\endcases      \eqno(\eqDeltaL)
$$


\subhead \sectCCG. Conjugacy classes in $GU(3,q^2)$ and $GL(3,q)$
\endsubhead
Recall that $GU(3,q^2)$
is the group of $3\times 3$ matrices $A$
with coefficients in the finite field $\GF{q^2}$ such that $A^*A=I$ where
$A^*=\overline{A^t}$ and $z\mapsto \bar z$ is the Frobenius automorphism
of $\GF{q^2}$ defined by $z\mapsto z^q$.

We set $\Omega=\{z\in\GF{q^2}\,|\,z^{q\pm1}=1\}$, i.~e., $\Omega$ is
the multiplicative group $\GF{q}^*$ when $G=GL$ and $\Omega$ is
``the unit circle''
$\UC=\{\,z\in\GF{q^2}\,|\,z\bar z=1\}$ when $G=GU$.

We fix a multiplicative generator $\tau$ of $\GF{q^6}^*$ and we
set $\rho = \tau^{q^4+q^2+1}$ (a generator of $\GF{q^2}^*$),
$\omega = \rho^{q\mp1}$ (a generator of $\Omega$),
and $\theta=\tau^{q^3\mp1}$.

The conjugacy classes in $GL(n,q)$ are determined by the Jordan normal form (JNF).
The conjugacy classes in $GU(n,q^2)$ have been computed in [\refEnnolaConj] and [\refWallG].
Each conjugacy class of $GU(n,q^2)$ is the intersection of $GU(n,q^2)$ with
a conjugacy class of $GL(n,q^2)$, so, it is determined by JNF.
The classes of $GL$ and those of $GU$
(represented by JNF in $GL(3,q^6)$) are
listed in Table \tabCC\ which, for the reader's convenience, we reproduce from [\refEnnola].
For an integer $k$, we denote the set $\{1,\dots,k\}$ by $[k]$.
We set $R_{q^2-1}=\{k\in[q^2-1]\,|\,k\not\equiv0\mod q\mp1\}$ and
$R_{q^3\pm1}=\{k\in[q^3\pm1]\,|\,k\not\equiv0\mod q^2\mp q+1\}$.

\midinsert
\noindent Table \tabCC. Conjugacy classes in $G$
\smallskip
\vbox{\offinterlineskip
\hrule
\halign{&\vrule#&\strut\;\;#\hfill\cr
height2pt&\omit&&\omit&&\omit&&\omit&&\omit&\cr
& \lower5pt\hbox{Class} &
& \lower5pt\hbox{JNF over $\GF{q^6}$}&
& \lower5pt\hbox{$\det$} &
& \lower5pt\hbox{class size} &
& range of the&
\cr
&\omit&&\omit&&\omit&&\omit&
&parameters&
\cr
height2pt&\omit&&\omit&&\omit&&\omit&&\omit&\cr
\noalign{\hrule}
height2pt&\omit&&\omit&&\omit&&\omit&&\omit&\cr
& $C_1^{(k)}$ &
& $\left(\smallmatrix\omega^k&0&0\\0&\omega^k&0\\0&0&\omega^k\endsmallmatrix\right)$ &
& $\omega^{3k}$ &
& $1$ &
& $k\in[q\pm1]$ &
\cr
height2pt&\omit&&\omit&&\omit&&\omit&&\omit&\cr
\noalign{\hrule}
height2pt&\omit&&\omit&&\omit&&\omit&&\omit&\cr
& $C_2^{(k)}$ &
& $\left(\smallmatrix\omega^k&0&0\\1&\omega^k&0\\0&0&\omega^k\endsmallmatrix\right)$ &
& $\omega^{3k}$ &
& $(q\mp1)(q^3\pm1)$ &
& $k\in[q\pm1]$ &
\cr
height2pt&\omit&&\omit&&\omit&&\omit&&\omit&\cr
\noalign{\hrule}
height3pt&\omit&&\omit&&\omit&&\omit&&\omit&\cr
& $C_3^{(k)}$ &
& $\left(\smallmatrix\omega^k&0&0\\1&\omega^k&0\\0&1&\omega^k\endsmallmatrix\right)$ &
& $\omega^{3k}$ &
& $q(q^2-1)(q^3\pm1)$ &
& $k\in[q\pm1]$ &
\cr
height2pt&\omit&&\omit&&\omit&&\omit&&\omit&\cr
\noalign{\hrule}
height2pt&\omit&&\omit&&\omit&&\omit&&\omit&\cr
& $C_4^{(k,l)}$ &
& $\left(\smallmatrix\omega^k&0&0\\0&\omega^k&0\\0&0&\omega^l\endsmallmatrix\right)$ &
& $\omega^{2k+l}$ &
& $q^2(q^2\mp q+1)$ &
& $(k,l)\in[q\pm1]^2$, $k\ne l$ &
\cr
height2pt&\omit&&\omit&&\omit&&\omit&&\omit&\cr
\noalign{\hrule}
height3pt&\omit&&\omit&&\omit&&\omit&&\omit&\cr
& $C_5^{(k,l)}$ &
& $\left(\smallmatrix\omega^k&0&0\\1&\omega^k&0\\0&0&\omega^l\endsmallmatrix\right)$ &
& $\omega^{2k+l}$ &
& $q^2(q\mp1)(q^3\pm1)$ &
& $(k,l)\in[q\pm1]^2$, $k\ne l$ &
\cr
height2pt&\omit&&\omit&&\omit&&\omit&&\omit&\cr
\noalign{\hrule}
height2pt&\omit&&\omit&&\omit&&\omit&&\omit&\cr
& $\!C_6^{(k,l,m)}\!$ &
& $\left(\smallmatrix\omega^k&0&0\\0&\omega^l&0\\0&0&\omega^m\endsmallmatrix\right)$ &
& $\!\omega^{k+l+m}\!$ &
& $\!q^3(q\mp1)(q^2\mp q+1)$\; &
& $1\le k<l<m\le q\pm1$ &
\cr
height2pt&\omit&&\omit&&\omit&&\omit&&\omit&\cr
\noalign{\hrule}
height3pt&\omit&&\omit&&\omit&&\omit&&\omit&\cr
& $C_7^{(k,l)}$ &
& $\left(\smallmatrix\omega^k&0&0\\0&\rho^l&0\\0&0&\rho^{\mp ql}\endsmallmatrix\right)$ &
& $\omega^{k\mp l}$ &
& $q^3(q^3\mp1)$ &
& $\matrix (k,l)\in[q\pm1]\times R_{q^2-1}\\ C_7^{(k,l)}=C_7^{(k,\mp ql)}\endmatrix$ &
\cr
height3pt&\omit&&\omit&&\omit&&\omit&&\omit&\cr
\noalign{\hrule}
height3pt&\omit&&\omit&&\omit&&\omit&&\omit&\cr
& $C_8^{(k)}$ &
& $\!\!\left(\smallmatrix\theta^{k}&0&0\\0&\theta^{q^2k}&0\\0&0&\theta^{q^4k}\endsmallmatrix\right)$ &
& $\omega^{k}$ &
& $q^3(q\pm1)^2(q\mp1)$ &
& $\matrix k\in R_{q^3\pm1}\\ C_8^{(k)}=C_8^{(q^2k)}=C_8^{(q^4k)}\endmatrix$ &
\cr
height3pt&\omit&&\omit&&\omit&&\omit&&\omit&\cr
}\hrule}
\endinsert


\subhead \sectCCS. Conjugacy classes in $SU(3,q^2)$ and $SL(3,q)$
\endsubhead
If $3$ does not divide $q\pm1$, then $G=S\times Z(G)$ where $Z(G)=C_1\cong\UC$ is the
center of $G$, and hence, the classes of $S$ are just those classes
of $G$ which are contained in $S$.

Let $q=3r\mp1$.
In this case, the splitting of conjugacy classes in $SL$ is described in
[\refDickson; Ch.~11, \S224] (see also [\refWallSL]). As stated in [\refSF], ``it can be shown
that the same splitting takes place in the unitary case''.
Each of $C_3^{(k)}$, $k=0,r,2r$, splits into three classes which we denote by
$C_3^{(k,l)}$, $l=0,1,2$. The class $C_3^{(k,l)}$ in $SU(3,q^2)$ (resp. in $SL(3,q)$)
consists of matrices which are conjugate in $SL(3,q^2)$ (resp. in $SL(3,q)$)
to\footnote{there is a misprint here in [\refSF].}
$$
   \left(\matrix \omega^k & 0        & 0       \\
                 z^l      & \omega^k & 0       \\
                 0        & 1        & \omega^k
   \endmatrix\right),
     \qquad  z = \cases \rho, & S=SU(3,q^2),\\ \omega, & S=SL(3,q).\endcases
$$
Other conjugacy classes of $G$ contained in $S$ are conjugacy
classes of $S$.

\proclaim{ Proposition \propSUinv }
If $A\in C_3^{(k,l)}$, then $A^{-1}\in C_3^{(-k,l)}$ and $\omega^{k'}A\in C_3^{(k+k',l)}$.
\endproclaim

\remark{ Remark } Each conjugacy class of $SU$ is the intersection of $SU$ with a conjugacy
class of $SL(3,q^2)$. The situation is quite different for $SU(2,q^2)$, see \S\sectNTwoS.
\endremark


\subhead \sectNotation. Notation for eigenvalues
\endsubhead
We denote the union of the conjugacy classes
$C_i^{(\dots)}$ by $C_i$, $i=1,\dots,8$.
We denote the number of distinct eigenvalues of matrices from $C_i$ by $n_i$
and the number of distinct eigenvalues belonging to $\Omega$
by $n'_i$. So, we have
$$
\split
&n_1=n_2=n_3=1,\qquad n_4=n_5=2,\qquad n_6=n_7=n_8=3;\\
&n'_i=n_i\quad (i=1,\dots,6);\qquad n'_7=1,\qquad n'_8=0.
\endsplit
$$

We denote the multiplicity of an eigenvalue $\lambda$ of a matrix $A$ by $m_A(\lambda)$.
Let $A\in C_i$.
We denote the eigenvalues of $A$ by
$\lambda_1=\lambda_1(A),\dots,\lambda_{n_i}=\lambda_{n_i}(A)$. We number them
so that 
$$
   m_A(\lambda_1)\ge\dots\ge m_A(\lambda_{n_i})\quad\text{and}\quad
   \lambda_1,\dots,\lambda_{n'_i}\in\UC.                             \eqno(\eqEV)
$$

For an $m$-tuple of matrices
$\vec A=(A_1,\dots,A_m)$, $A_\nu\in C_{i_\nu}$, $\nu=1,\dots,m$,
we use the multi-index notation:
$$
   \vec a=(a_1,\dots,a_m),\quad
   [\vec{n}]=[n_{i_1}]\times\dots\times[n_{i_m}],\quad
   [\vec{n'}]=[n'_{i_1}]\times\dots\times[n'_{i_m}],\quad
$$
(recall that $[k]$ stands for $\{1,\dots,k\}$) and for $\vec a\in[\vec n]$ we set
$$
   \lambda_{\vec a} = \lambda_{a_1}(A_1)\dots\lambda_{a_m}(A_m),\qquad
   \delta_{\vec a} = \delta_{\vec a}(\vec A)=\cases 1, & \lambda_{\vec a}=1,\\
                            0, & \lambda_{\vec a}\ne 1. \endcases
$$

In this notation, the rank condition (\eqRk) for $\vec A=(A_1,A_2,A_3)$ takes the form
$$
        \sum_{a=1}^{n'_{i_3}} \delta_{1,1,a} > 0
\qquad
	\text{if $\{i_1,i_2\}\subset\{2,4\}$}
						             \eqno(\eqRk')
$$


\subhead \sectMainTh. Statement of main results
\endsubhead
In Theorems \thGU\ and \thSU, we restrict ourselves by the case when $A_1,\dots,A_m$
are non-scalar and $m\ge3$. To reduce the general case to this one, it is enough to
know the class of the inverse of a given matrix and the class of its multiple by a scalar.
For $G$, this is clear from JNF; for $S$, the answer is given in Proposition \propSUinv\ in \S\sectCCS.

\proclaim{ Theorem \thGU } 
Let $A_1,\dots,A_m\in G\setminus C_1$, $m\ge 3$,  satisfy {\rm(\eqDet)}
and {\rm(\eqRk)}. Let $A_\nu\in C_{i_\nu}$, $\nu=1,\dots,m$.

\medskip
{\rm(a)}.
If $G=GU$, we suppose that one of the following conditions $(i)$--$(vii)$ holds:
\roster
\item"$(i)$"
        $m=3$, $i_1\in\{6,7\}$, $i_2\in\{3,5\}$, $i_3\in\{2,4\}$, and $\delta_{111}=1$;
\item"$(ii)$"
        $m=3$, $i_1=5$, $i_2\in\{3,5\}$, $i_3\in\{2,4\}$, and $\delta_{211}=1$;
\item"$(iii)$"
        $m=3$, $i_1=i_2\in\{6,8\}$, $i_3=2$, and $\delta_{111}\delta_{221}\delta_{331}=1$
        {\rm(when $i_1=i_2=8$, the last condition is equivalent to $\delta_{111}=1$)};
\item"$(iv)$"
        $m=3$, $(i_1,i_2,i_3)=(3,2,2)$ or $(4,4,2)$;
\item"$(v)$"
        $m=3$, $(i_1,i_2,i_3)=(5,4,4)$, and $\delta_{112}\delta_{121}\delta_{211}=1$
        {\rm(}see Remark \remGUbis{\rm)};
\item"$(vi)$"
        $m=4$, $(i_1,i_2,i_3,i_4)=(3,2,2,2)$, and $\delta_{1111}=1$;
\item"$(vii)$"
        $m=4$, $(i_1,i_2,i_3,i_4)=(4,4,4,2)$, and $\delta_{1121}\delta_{1211}\delta_{2111}=1$
        {\rm(}see Remark \remGUbis{\rm)}.
\endroster
If $G=GL$, we suppose that one of the following conditions $(viii)$--$(ix)$ holds:
\roster
\item"$(viii)$"
        $m=3$, $(i_1,i_2,i_3)=(8,8,2)$ and $\delta_{111}=1$;
\item"$(ix)$"
        $q=2$, $m=3$, $(i_1,i_2,i_3)=(8,8,3)$ and $A_1^G=A_2^G$.
\endroster
Then $I\not\in A_1^G\dots A_m^G$.

\medskip
{\rm(b)}
Suppose that none of the conditions of Part {\rm(a)} holds for
any permutation of $A_1,\dots,A_m$ and for any renumbering of the
eigenvalues of the matrices under the restrictions {\rm(\eqEV)}.
In the case $G=GU$, we suppose also that $q\ne2$.
Then $I\in A_1^G\dots A_m^G$.
\endproclaim

\remark{ Remark \remGU }
In Table \tabExcept\ we present the list of all the cases when (\eqDet) can be satisfied for
non-constant matrices $A_1,\dots,A_m\in G$, $m\ge 3$,
$A_i\in C_{i_\nu}$, $i_\nu\ge2$, but
$I\not\in A_1^G\dots A_m^G$ for $q>2$.
The cases marked by asterisk concern the both groups $GU$ and $GL$;
as stated in \S\sectDetRk, in all of them except $(8,8,2)$ the rank condition (\eqRk)
is not satisfied. The cases not marked by asterisk concern only $GU$.
\endremark

\remark{ Remark \remGUbis }
Conditions $(v)$ and $(vii)$ in Theorem \thGU\ mean that
$(\mu_1A_1^G,\dots,\mu_m A_m^G)$ is
$\big( C_5^{(-k,0)}, C_4^{(0,k)}, C_4^{(0,k)}\big)$ or
$\big( C_4^{(0,k)}, C_4^{(0,k)}, C_4^{(0,k)}, C_2^{(-k)}\big)$
for some $k\in\{1,\dots,q\}$ and for some $\mu_i\in\Omega$ with $\mu_1\dots\mu_m=1$.
\endremark

\midinsert
\noindent Table \tabExcept. Cases when $\det A_1\dots A_m=1$,
$I\not\in A_1^G\dots A_m^G$ for $q>2$
(see Remark \remGU)
\smallskip
\vbox{\offinterlineskip
\hrule
\halign{&\vrule#&\strut\quad#\hfill\cr
height3pt&\omit&\omit&\omit&&\omit&\omit&\omit&\cr
& $\!\!(i_1,\dots,i_m)$ &\omit&\omit&
& $\!\!(i_1,\dots,i_m)$ &\omit&\omit&\cr
height3pt&\omit&\omit&\omit&&\omit&\omit&\omit&\cr
\noalign{\hrule}
height3pt&\omit&\omit&\omit&&\omit&\omit&\omit&\cr
& $(2,2,2)$*&\omit& $\delta_{111}=0$ &
& $(6,4,2)$*&\omit& $\delta_{111}+\delta_{211}+\delta_{311}=0$ &
\cr
& $(3,2,2)$*&\omit& $\delta_{111}=0$ &
& $(6,4,3)$ &\omit& $\delta_{111}+\delta_{211}+\delta_{311}=1$ &
\cr
& $(3,2,2)$ &\omit& $\delta_{111}=1$ &
& $(6,4,4)$*&\omit& $\delta_{111}+\delta_{211}+\delta_{311}=0$ &
\cr
& $(4,2,2)$*&\omit& &
& $(6,5,2)$ &\omit& $\delta_{111}+\delta_{211}+\delta_{311}=1$ &
\cr
& $(4,3,2)$*&\omit& &
& $(6,5,4)$ &\omit& $\delta_{111}+\delta_{211}+\delta_{311}=1$ &
\cr
& $(4,4,2)$*&\omit& $\delta_{111}=0$ &
& $(6,6,2)$ &\omit& $\sum_{\alpha\in S_3}\delta_{11^\alpha1}\delta_{22^\alpha1}\delta_{33^\alpha1}=1$ &
\cr
& $(4,4,2)$ &\omit& $\delta_{111}=1$ &
& $(7,2,2)$*&\omit& $\delta_{111}=0$ &
\cr
& $(4,4,3)$*&\omit& $\delta_{111}=0$ &
& $(7,3,2)$ &\omit& $\delta_{111}=1$ &
\cr
& $(4,4,4)$*&\omit& $\delta_{111}+\delta_{112}\delta_{121}\delta_{211}=0$ &
& $(7,4,2)$*&\omit& $\delta_{111}=0$ &
\cr
& $(5,2,2)$*&\omit& $\delta_{211}=0$ &
& $(7,4,3)$ &\omit& $\delta_{111}=1$ &
\cr
& $(5,3,2)$ &\omit& $\delta_{211}=1$ &
& $(7,4,4)$*&\omit& $\delta_{111}=0$ &
\cr
& $(5,4,2)$*&\omit& $\delta_{111}+\delta_{211}=0$ &
& $(7,5,2)$ &\omit& $\delta_{111}=1$ &
\cr
& $(5,4,3)$ &\omit& $\delta_{211}=1$ &
& $(7,5,4)$ &\omit& $\delta_{111}=1$ &
\cr
& $(5,4,4)$*&\omit& $\delta_{111}+\delta_{211}=0$ &
& $(8,2,2)$*&\omit&  &
\cr
& $(5,4,4)$ &\omit& $\delta_{211}\delta_{121}\delta_{112}=1$ &
& $(8,4,2)$*&\omit&  &
\cr
& $(5,5,2)$ &\omit& $\delta_{211}=1$ &
& $(8,4,4)$*&\omit&  &
\cr
& $(5,5,4)$ &\omit& $\delta_{211}=1$ &
& $(8,8,2)$*&\omit& $\delta_{111}+\delta_{121}+\delta_{131}=1$ &
\cr
& $(6,2,2)$*&\omit& $\delta_{111}+\delta_{211}+\delta_{311}=0$ &
& $(3,2,2,2)$ &\omit& $\delta_{1111}=1$ &
\cr
& $(6,3,2)$ &\omit& $\delta_{111}+\delta_{211}+\delta_{311}=1$ &
& $(4,4,4,2)$ &\omit& $\delta_{1121}\delta_{1211}\delta_{2111}=1$ &
\cr
height3pt&\omit&\omit&\omit&&\omit&\omit&\omit&\cr
}\hrule}
\endinsert

The case of $q=2$ also is treated completely in
Propositions \propSmallqTwo\ and
\propSmallqThree\ (for $GU$) and in Corollary
\corSmallqSU\ (for $SU$). 

If $3$ does not divide $q\pm1$, then $G\cong S\times\Omega$, thus
the class product problem in $S$ reduces to that in $G$.
Otherwise (when $3|q\pm1$) the solution
is as follows.

\proclaim{ Theorem \thSU }
Let $A_1,\dots,A_m$, $m\ge 3$, be as in Theorem \thGU.
We suppose in addition that $q=3r\mp1$ and $A_1,\dots,A_m\in S$,
recall that $S$ is $SU(3,q^2)$ or $SL(3,q)$.

\medskip
{\rm(a)}. Suppose that $m=3$.

\noindent
If $S=SU$, we suppose that
\roster
\item"$(i)$"
              $i_1=i_2=3$, $i_3\in\{2,4\}$,
              $A_1\in C_3^{(k_1,l_1)}$, $A_2\in C_3^{(k_2,l_2)}$,
               $l_1\ne l_2$.
\endroster
If $S=SL$, we suppose that one of the following conditions $(ii)$--$(v)$ holds:
\roster
\item"$(ii)$"
              $(i_1,i_2,i_3)=(3,3,2)$,
              $A_1\in C_3^{(k_1,l_1)}$, $A_2\in C_3^{(k_2,l_2)}$,
              $l_1\ne l_2$, and $\delta_{111}=0$;
\item"$(iii)$"
              $(i_1,i_2,i_3)=(3,3,4)$,
              $A_1\in C_3^{(k_1,l_1)}$, $A_2\in C_3^{(k_2,l_2)}$,
              $l_1\ne l_2$;
\item"$(iv)$"
              $q=4$, $(i_1,i_2,i_3)=(3,3,3)$,
              $A_\nu\in C_3^{(k_\nu,l_\nu)}$, $\nu=1,2,3$,
              $l_1=l_2\ne l_3$, and $\delta_{111}=1$;
\item"$(v)$"
              $q=4$, $(i_1,i_2,i_3)=(3,3,3)$,
              $A_\nu\in C_3^{(k_\nu,l_\nu)}$, $\nu=1,2,3$,
              $l_1=l_2=l_3$, and $\delta_{111}=0$.
\endroster
Then $I\not\in A_1^S A_2^S A_3^S$.

\medskip
{\rm(b)}. Suppose that $q>2$ and $I\in A_1^G\dots A_m^G$. Suppose that
for any permutation of $A_1,\dots,A_m$,
the hypothesis of Part {\rm(a)} is not satisfied.
Then $I\in A_1^S\dots A_m^S$.
\endproclaim

If $3$ does not divide $q\pm1$, then $PG=PS=S$. If $3$ divides $q\pm1$, the
solution of the class product problem for $PG$ and $PS$ is as follows.
Let $\tilde C_i^{(\dots)}$ be the conjugacy class of $PG$ or $PS$
corresponding to $C_i^{(\dots)}$.

\proclaim{ Corollary \corPGU }
Let $q=3r\mp 1$, $q\ne2$. If $m\ge4$ {\rm(}resp. $m\ge3${\rm)},
then the product of any $m$-tuple of nontrivial conjugacy classes
of $PGU$ {\rm(}resp. $PGL${\rm)} contains the identity matrix.
All triples of nontrivial conjugacy classes of $PGU$ which have representatives
in $GU$ satisfying {\rm(\eqDet)} and {\rm(\eqRk)}, but
whose product does not contain the identity matrix, are
$$
\def\tC{\tilde C}
\xalignat5
&(i)&&\tC_3^{(0)} &&\tC_2^{(0)} &&\tC_2^{(0)} &\qquad\qquad& \\
&(ii)&&\tC_2^{(0)} &&\tC_4^{(0,k)} &&\tC_4^{(0,-k)}& &k=1,\dots q;\\
&(iii)&&\tC_5^{(0,k)} &&\tC_4^{(0,k)} &&\tC_4^{(0,k)}& &k=1,\dots,q,\;\;k\not\in\{r,2r\};\\
&(iv)&&\tC_6^{(0,r,2r)}&&\tC_3^{(0)}   &&\tC_2^{(0)}    & &\\
&(v)&&\tC_6^{(0,r,2r)}&&\tC_5^{(0,k)} &&\tC_4^{(0,-k)}& &k=1,\dots,q;\\
&(vi)&&\tC_6^{(0,r,2r)}&&\tC_6^{(0,r,2r)}&&\tC_2^{(0)}    &
\endxalignat
$$
\endproclaim

\proclaim{ Corollary \corPSU }
Let $q=3r\mp 1$, $q\ne2$. If $m\ge4$, then the product of any $m$-tuple of nontrivial conjugacy classes
of $PS$ contains the identity matrix.
All triples of nontrivial conjugacy classes which have representatives
in $S$ satisfying {\rm(\eqRk)}, but
whose product does not contain the identity matrix, are 
$$
\def\tC{\tilde C}
\xalignat5
&(i)&&\tC_3^{(0,l)} &&\tC_2^{(0)} &&\tC_2^{(0)} & &l=0,1,2;\\
&(ii)&&\tC_2^{(0)} &&\tC_4^{(k,-2k)} &&\tC_4^{(-k,2k)}& &k=1,\dots r-1;\\
&(iii)&&\tC_5^{(k,-2k)} &&\tC_4^{(k,-2k)} &&\tC_4^{(k,-2k)}& &k=1,\dots,r-1,\;\;3k\not\in\{r,2r\};\\
&(iv)&&\tC_6^{(0,r,2r)}&&\tC_3^{(0,l)}   &&\tC_2^{(0)}    & &l=0,1,2;\\
&(v)&&\tC_6^{(0,r,2r)}&&\tC_5^{(k,-2k)} &&\tC_4^{(k,-2k)}& &k=1,\dots,r-1;\\
&(vi)&&\tC_6^{(0,r,2r)}&&\tC_6^{(0,r,2r)}&&\tC_2^{(0)}    &
\endxalignat
$$
$$
\def\tC{\tilde C}
\xalignat5
&(vii)&&\tC_3^{(0,l_1)} &&\tC_3^{(0,l_2)} &&\tC_2^{(0)}    & \qquad\qquad&0\le l_1<l_2\le 2;\\
&(viii)&&\tC_3^{(0,l_1)} &&\tC_3^{(0,l_2)} &&\tC_4^{(k,-2k)}& &0\le l_1<l_2\le 2,\;\; k=1,\dots,r-1.
\endxalignat
$$
in the case $PS=PSU$, and only the triples $(viii)$ in the case $PS=PSL$. 
\endproclaim


\subhead \sectCN. Covering number and extended covering number
\endsubhead
Let $\Gamma$ be a group. The {\it covering number} of $\Gamma$ is the minimal integer $m$
such that for any nontrivial conjugacy class $c$, we have $c^m=\Gamma$. It is denoted by $\cn(\Gamma)$.
The {\it extended covering number} of $\Gamma$ is the minimal integer $m$
such that for any nontrivial conjugacy classes $c_1,\dots,c_m$ we have $c_1\dots c_m=\Gamma$.
Covering numbers were studied in [\refAH, \refLSh].

\proclaim{ Corollary \corECN }

$\cn(PSL)=3$ and $\ecn(PSL)=4$;

$\cn(PSU)=3$ and $\ecn(PSU)=4$ if $\gcd(q+1,3)=3$ and $q\ne2$;

$\cn(PSU)=4$ and $\ecn(PSU)=5$ if $\gcd(q+1,3)=1$.
\smallskip
\endproclaim

\remark{ Remark \remKarniLev}
Karni [\refKarni] computed
the numbers $\cn(PS)$ and $\ecn(PS)$ for $q=3,4,5$;
Lev [\refLev] proved that $\cn(PSL(n,K))=n$ for any $n\ge 3$ and for any field $K$ which
has more than $3$ elements.
\endremark


\head \sectG.  Class products in $GU(3,q^2)$ and $GL(3,q)$. Proof of Theorem \thGU
\endhead

\subhead \sectCT. The character tables of $GU(3,q^2)$ and $GL(3,q)$
\endsubhead
%
In this section we represent the character table of $G$ (see [\refEnnola])
in a form convenient to apply (\eqCSF).
The irreducible characters of $G$ divide into $8$ series parametrized
by the same sets of parameters as the conjugacy classes.
We denote the dimension of the irreducible representations corresponding
to the $j$-th series by $d_j$. So,
$$
\xalignat4
d_1&=1        & d_3&=q^3        & d_5&=q(q^2\mp q+1)     & d_7&=q^3\pm1 \\
d_2&=q^2\mp q & d_4&=q^2\mp q+1 & d_6&=(q\mp1)(q^2\mp q+1) & d_8&=(q\pm1)(q^2-1)
\endxalignat
$$
The characters $\chi_{d_i}^{(t,\dots)}$, $i=1,\dots,8$,
are irreducible and pairwise distinct only for some values of the parameters
$t$, $u$, $v$, but we define them by the same formulas for any values of the
parameters. Recall that for an integer $n$, we denote the set $\{1,\dots,n\}$ by $[n]$.
Let 
$$
\xalignat 2
X_j &=\{\chi_{d_j}^{(t)}    \,|\, t     \in[q\pm1]  \},\;j=1,2,3,&
X'_j&=\{\chi_{d_j}^{(t,t)}  \,|\, t     \in[q\pm1]  \},\;j=4,5,
\\
X_j &= \{\chi_{d_j}^{(t,u)}  \,|\,(t,u)  \in[q\pm1]^2 \},\;j=4,5,&
X'_6&= \{\chi_{d_6}^{(t,u,u)}\,|\,(t,u)  \in[q\pm1]^2 \},
\\
X_6 &= \{\chi_{d_6}^{(t,u,v)}  \,|\,(t,u,v)\in[q\pm1]^3 \}.&
X'_7&= \{\chi_{d_7}^{(t,(1\mp q)u)}\,|\,  (t,u)\in[q\pm1]   \},
\\
X_7 &= \{\chi_{d_7}^{(t,u)}      \,|\,(t,u)  \in[q\pm1]\times[q^2-1]\},&
X'_8&= \{\chi_{d_8}^{((q^2\mp q+1)t)}\,|\,   t   \in[q\pm1]  \},
\\
X_8  &= \{\chi_{d_8}^{(t)}    \,|\,t\in[q^3\pm1] \},&
X''_6&= \{\chi_{d_6}^{(t,t,t)}\,|\,t\in[q\pm1]   \}
\endxalignat
$$
and
$\Xi_1 = \{X_1,X_2,X_3,X'_4,X'_5,X''_6,X'_8\}$,
$\Xi_2 = \{X_4,X_5,X'_6,X'_7\}$,
$\Xi_3 = \{X_6\}$,
$\Xi_4 = \{X_7\}$,
$\Xi_5 = \{X_8\}$,
$\Xi=\Xi_1\cup \dots\cup \Xi_5$.
It is clear that if $E$ is any expression depending on a character of $G$, then
$$
    \sum_{\chi\in \Irr(G)} E(\chi)
   = \sum_{X\in\Xi} s(X) \sum_{\chi\in X}E(\chi) \eqno(\eqSumOverIrr)
$$
where the symmetry factors $s(X)$ are given in Tables \tabCT.1 and \tabCT.2.

We fix a homomorphism of multiplicative groups $f:\GF{q^6}^*\to\Bbb C^*$
which takes $\tau$ to $\exp(2\pi i/(q^6-1))$, thus,
$$
   f(\omega)=e^{2\pi i/(q\pm1)},\quad
   f(\rho)=e^{2\pi i/(q^2-1)},\quad
   f(\theta)=e^{2\pi i/(q^3\pm1)}.
$$

Let $A\in C_i$ and let $\lambda_1,\dots,\lambda_{n_i}$ be its eigenvalues
numbered as in (\eqEV).
Then
$$
 \chi^{(t)}(A) = c_i^X f(\det A)^t,\qquad \chi^{(t)}\in X\in\Xi_1
$$
$$
 \chi^{(t,u)}(A) = \sum_{a=1}^{n'_i} c_{i,a}^X
   f(\lambda_a)^t f(\lambda_a^{-1}\det A)^u,
 \qquad \chi^{(t,u)}\in X\in\Xi_2,
$$
$$
     \chi_{d_6}^{(t,u,v)}(A) = \sum_{\alpha\in\Cal A_{6,i}} c^{X_6}_i
          f(\lambda_{\alpha(1)}^t
            \lambda_{\alpha(2)}^u
            \lambda_{\alpha(3)}^v),
$$
$$
    \chi_{d_7}^{(t,u)}(A) = \sum_{\alpha\in\Cal A_{7,i}} c^{X_7}_i
          f(\lambda_{\alpha(1)}^t \lambda_{\alpha(2)}^u),
\qquad\qquad
    \chi_{d_8}^{(t)}(A) = \sum_{a=1}^{n_i} c^{X_8}_i f(\lambda_a^t).
$$
where $\Cal A_{6,i}$ and $\Cal A_{7,i}$ are sets of
triples $\alpha=(\alpha(1),\alpha(2),\alpha(3))$ and
pairs $(\alpha(1),\alpha(2))$ respectively defined by
$$
\xalignat3
&\Cal A_{6,i} = \{(1,1,1)\},&&\Cal A_{7,i}=\{(1,1)\},&& i=1,2,3,\\
&\Cal A_{6,i} = \{(2,1,1),(1,2,1),(1,1,2)\},&&\Cal A_{7,i}=\{(2,1)\},&&i=4,5,\\
&\Cal A_{6,6} = S_3, && \Cal A_{7,7}= \{(1,2),(1,3)\},\\
&\Cal A_{6,7}=\Cal A_{6,8} = \varnothing, && \Cal A_{7,6}=\Cal A_{7,8}=\varnothing.
\endxalignat
$$
The coefficients $c_i^X$ and $c_{i,a}^X$ (the latter denoted just by
$c_i^X$ in the cases when $n'_i=1$)
are given in the Tables \tabCT.1 and \tabCT.2.

\midinsert
\noindent Table \tabCT.1
\smallskip
\vbox{\offinterlineskip
\hrule
\halign{&\vrule#&\strut\hfill\;\;\;#\;\;\hfill\cr
height3pt%
&\omit&
&\omit&\omit&\omit&\omit&\omit&\omit&\omit&\omit
&\omit&\omit&\omit&\omit&\omit&\omit&\omit&&\omit&\cr
& $X$ &
& $c_1^X$ &\omit& $c_2^X$ &\omit& $c_3^X$ &\omit& $c_4^X$ &\omit
& $c_5^X$ &\omit& $c_6^X$ &\omit& $c_7^X$ &\omit& $c_8^X$ && $s(X)$ &\cr
height3pt&\omit&
&\omit&\omit&\omit&\omit&\omit&\omit&\omit&\omit
&\omit&\omit&\omit&\omit&\omit&\omit&\omit&&\omit&\cr
\noalign{\hrule}
height3pt&\omit&
&\omit&\omit&\omit&\omit&\omit&\omit&\omit&\omit
&\omit&\omit&\omit&\omit&\omit&\omit&\omit&&\omit&\cr
& $X_1$ &
& $1$   &\omit& $1$ &\omit& $1$ &\omit& $1$ &\omit
& $1$   &\omit& $1$ &\omit& $1$ &\omit& $1$ && $1$ &
\cr
& $X_2$ &
& $d_2$ &\omit& $\mp q$ &\omit&  $0$ &\omit&  $1\mp q$ &\omit
&  $1$  &\omit&     $2$ &\omit&  $0$ &\omit&      $-1$ && $1$ &     
\cr
& $X_3$ &
& $d_3$ &\omit&    $0$ &\omit&    $0$ &\omit&    $q$ &\omit
&   $0$ &\omit& $\mp1$ &\omit& $\pm1$ &\omit& $\mp1$ && $1$ &     
\cr
& $X'_4$ &
& $d_4$ &\omit& $1\mp q$ &\omit&  $1$ &\omit& $2\mp q$ &\omit
&   $2$ &\omit&      $3$ &\omit&  $1$ &\omit&      $0$ && $-1$ &    
\cr
& $X'_5$ &
& $d_5$ &\omit&    $q$ &\omit&    $0$ &\omit& $2q\mp1$ &\omit
&  $-1$ &\omit& $\mp3$ &\omit& $\pm1$ &\omit&      $0$ && $-1$ &    
\cr
& $X''_6$ &
&   $d_6$ &\omit& $2q\mp1$ &\omit& $\mp1$ &\omit& $3q\mp3$ &\omit
&  $\mp3$ &\omit&   $\mp6$ &\omit&    $0$ &\omit&      $0$ && $1/3$ &   
\cr
& $X'_8$ &
& $d_8$ &\omit& $-q\mp1$ &\omit& $\mp1$ &\omit&    $0$ &\omit
&   $0$ &\omit&      $0$ &\omit&    $0$ &\omit& $\mp3$ && $-1/3$ &  
\cr
height3pt&\omit&
&\omit&\omit&\omit&\omit&\omit&\omit&\omit&\omit
&\omit&\omit&\omit&\omit&\omit&\omit&\omit&&\omit&\cr
\noalign{\hrule}
height3pt&\omit&
&\omit&\omit&\omit&\omit&\omit&\omit&\omit&\omit
&\omit&\omit&\omit&\omit&\omit&\omit&\omit&&\omit&\cr
& $X_6$ &
&  $d_6$ &\omit& $2q\mp1$ &\omit& $\mp1$ &\omit& $q\mp1$ &\omit
& $\mp1$ &\omit&   $\mp1$ &\omit&    $0$ &\omit&     $0$ && $1/6$ &
\cr
& $X_7$ &
& $d_7$ &\omit&    $\pm1$ &\omit&  $\pm1$ &\omit& $q\pm1$ &\omit
& $\pm1$ &\omit&      $0$ &\omit&  $\pm1$ &\omit&     $0$ && $1/2$ &
\cr
& $X_8$ &
& $d_8$ &\omit& $-q\mp1$ &\omit& $\mp1$ &\omit&    $0$ &\omit
&   $0$ &\omit&      $0$ &\omit&    $0$ &\omit& $\mp1$ && $1/3$ &
\cr
height3pt&\omit&
&\omit&\omit&\omit&\omit&\omit&\omit&\omit&\omit
&\omit&\omit&\omit&\omit&\omit&\omit&\omit&&\omit&\cr
}\hrule}
\endinsert

\midinsert
\noindent Table \tabCT.2
\smallskip
\vbox{\offinterlineskip
\hrule
\halign{&\vrule#&\strut\hfill\;\;#\;\;\hfill\cr
height3pt%
&\omit&
&\omit&\omit&\omit&\omit&\omit&\omit&\omit&\omit&\omit
&\omit&\omit&\omit&\omit&\omit&\omit&\omit&\omit&&\omit&\cr
& $X$ &
& $c_2^X$ &\omit& $c_2^X$ &\omit& $c_3^X$ &\omit& $c_{4,1}^X$ &\omit& $c_{4,2}^X$ &\omit
& $c_{5,1}^X$ &\omit& $c_{5,2}^X$ &\omit& $c_{6,a}^X$ &\omit& $c_7^X$ && $s(X)$ &\cr
height3pt&\omit&
&\omit&\omit&\omit&\omit&\omit&\omit&\omit&\omit&\omit
&\omit&\omit&\omit&\omit&\omit&\omit&\omit&\omit&&\omit&\cr
\noalign{\hrule}
height3pt&\omit&
&\omit&\omit&\omit&\omit&\omit&\omit&\omit&\omit&\omit
&\omit&\omit&\omit&\omit&\omit&\omit&\omit&\omit&&\omit&\cr
& $X_4$ &
& $d_4$ &\omit& $1\mp q$ &\omit& $1$ &\omit& $1\mp q$ &\omit& $1$ &\omit
&   $1$ &\omit&      $1$ &\omit& $1$ &\omit&      $1$ && $1$ &
\cr
& $X_5$ &
& $d_5$ &\omit& $q$ &\omit& $0$ &\omit& $q\mp1$ &\omit& $q$ &\omit
& $\mp1$ &\omit& $0$ &\omit& $\mp1$ &\omit& $\pm1$   && $1$ &
\cr
& $X'_6$ &
& $d_6$ &\omit& $2q\mp1$ &\omit& $\mp1$ &\omit& $2(q\mp1)$ &\omit& $q\mp1$ &\omit
& $\mp2$ &\omit& $\mp1$ &\omit& $\mp2$ &\omit& $0$ && $-1/2$ &
\cr
& $X'_7$ &
& $d_7$ &\omit& $\pm1$ &\omit& $\pm1$ &\omit& $0$ &\omit& $q\pm1$ &\omit
& $0$ &\omit& $\pm1$ &\omit& $0$ &\omit& $\pm2$ && $-1/2$ &
\cr
height3pt&\omit&
&\omit&\omit&\omit&\omit&\omit&\omit&\omit&\omit&\omit
&\omit&\omit&\omit&\omit&\omit&\omit&\omit&\omit&&\omit&\cr
}\hrule}
\endinsert


\subhead\sectSCF.  Structure constant formula for $GU(3,q^2)$ and $GL(3,q)$
\endsubhead
%
Let $A_1,\dots,A_m\in G$, $A_\nu\in C_{i_\nu}$, $\det A_1\dots A_m=1$.
We use the multi-index notation as explained in \S\sectNotation\ and we set also
$$
    \vec\Cal A_j =\Cal A_{j,i_1}\times\dots\times\Cal A_{j,i_m},\quad
    j=6,7.
$$
Substituting the formulas from \S\sectCT\ into (\eqSCF) and using
(\eqSumOverIrr),
we obtain
$$
   \bar N_G(A_1,\dots,A_m)
   =\Sigma_1+\dots+\Sigma_5
$$
where $\Sigma_i$ is the sum over $\Xi_i$:
$$
\split
  \Sigma_1 &= \sum_{X\in\Xi_1} s(X) \sum_{t=1}^{q\pm1}
      {c_{i_1}^X\dots c_{i_m}^X f\big(\det A_1\dots A_m\big)^t\over (c_{1}^X)^{m-2}}
   = (q\pm1)\sum_{X\in\Xi_1} {s(X)c_{i_1}^X\dots c_{i_m}^X\over (c_{1}^X)^{m-2}}
\endsplit
$$
$$
\split
   \Sigma_2 &= \sum_{X\in\Xi_2} s(X)\sum_{\vec a\in[\vec n']}
      {c_{i_1,a_1}^X\dots c_{i_m,a_m}^X\over (c_1^X)^{m-2}}
     \sum_{t=1}^{q\pm1}
       f\big( \lambda_{\vec a} \big)^t
         \sum_{u=1}^{q\pm1}
       f\big( \lambda_{\vec a}^{-1} \det A_1\dots A_m\big)^u
\\&
    =(q\pm1)^2\sum_{X\in\Xi_2} s(X)\sum_{\vec a\in[\vec n']}
      {c_{i_1,a_1}^X\dots c_{i_m,a_m}^X\delta_{\vec a}
               \over (c_1^X)^{m-2}}
\endsplit
$$
$$
\split
   \Sigma_3 &=
    {1\over6}\sum_{\vec\alpha\in\vec\Cal A_{6}}
        {c_{i_1}^{X_6}\dots c_{i_m}^{X_6}\over d_6^{m-2}}
        \sum_{t=1}^{q\pm1}f(\lambda_{\vec\alpha(1)})^t
        \sum_{u=1}^{q\pm1}f(\lambda_{\vec\alpha(2)})^u
        \sum_{v=1}^{q\pm1}f(\lambda_{\vec\alpha(3)})^v
\\&
     = {(q\pm1)^3\over6}\sum_{\vec\alpha\in\vec\Cal A_{6}}
        {c_{i_1}^{X_6}\dots c_{i_m}^{X_6}\over d_6^{m-2}}
         \delta_{\vec\alpha(1)}\delta_{\vec\alpha(2)}\delta_{\vec\alpha(3)}
\endsplit
$$
$$
\split
   \Sigma_4 &=
    {1\over2}\sum_{\vec\alpha\in\vec\Cal A_{7}}
        {c_{i_1}^{X_7}\dots c_{i_m}^{X_7}\over d_7^{m-2}}
        \sum_{t=1}^{q\pm1}f(\lambda_{\vec\alpha(1)})^t
        \sum_{u=1}^{q^2-1}f(\lambda_{\vec\alpha(2)})^u
\\&
     = {(q\pm1)(q^2-1)\over2}\sum_{\vec\alpha\in\vec\Cal A_{7}}
        {c_{i_1}^{X_7}\dots c_{i_m}^{X_7}\over d_7^{m-2}}
         \delta_{\vec\alpha(1)}\delta_{\vec\alpha(2)}
\endsplit
$$
$$
\split
   \Sigma_5 &=
    {1\over3}\sum_{\vec a\in[\vec n]}
        {c_{i_1}^{X_8}\dots c_{i_m}^{X_8}\over d_8^{m-2}}
        \sum_{t=1}^{q^3\pm1}f(\lambda_{\vec a})^t
     = {(q^3\pm1)\over3}\sum_{\vec a\in[\vec n]}
        {c_{i_1}^{X_8}\dots c_{i_m}^{X_8}\over d_8^{m-2}}
         \delta_{\vec a}
\endsplit
$$


\subhead\sectTP. Structure constants for triple
                 products in $GU(3,q^2)$ and $GL(3,q)$ 
\endsubhead
%
%
Using the formulas from \S\sectSCF, we computed
the structure constants for all triples $(i_1,i_2,i_3)$.
To write down the result in a compact form, we introduce
the following notation.
We define $\vec\Cal A^*_6$ as the quotient of $\vec\Cal A_6$
by the action of the symmetric group $S_3$ defined by
$\vec\alpha^\pi=(\alpha_1^\pi,\dots,\alpha_m^\pi)$
where $\alpha_\nu^\pi=(\alpha_\nu(1^\pi),
\alpha_\nu(2^\pi),\alpha_\nu(3^\pi))$.
Similarly, we define $\vec\Cal A^*_7$ as the quotient
of $\vec\Cal A_7$ by the action of $\Bbb Z_2$ which exchanges
the elements of $\Cal A_{7,7}$.
Given $\vec a\in[\vec n']$, let
$|\vec a|$ be the number of $\nu$ such that $a_\nu=1$ and
$i_\nu\in\{4,5\}$.
We set
$$
 \Delta = \sum_{\vec\alpha\in\vec\Cal A^*_6}
         \delta_{\vec\alpha(1)}
         \delta_{\vec\alpha(2)}
         \delta_{\vec\alpha(3)} + \sum_{\vec\alpha\in\vec\Cal A^*_7}
         \delta_{\vec\alpha(1)}\delta_{\vec\alpha(2)},
\qquad
    \Delta_a = \sum_
    {\vec a\in[\vec n'],|\vec a|=a}
    \delta_{\vec a}.
$$
We set also
$$
    \Delta'=\sum_{\vec a\in[\vec n]}\delta_{\vec a}.
$$

\def\evi{\text{ev.}}
\def\TT{tbl.~\tabTPExcept}
\def\STP{\S\sectTPtwo}
\midinsert
\noindent Table \tabTripleProd. Structure constants for $GU(3,q^2)$
$\left(\smallmatrix\delta_L=0\\\pm=+\endsmallmatrix\right)$
and $GL(3,q)$
$\left(\smallmatrix\delta_L=1\\\pm=-\endsmallmatrix\right)$
\smallskip
\vbox{\offinterlineskip
\hrule
\halign{&\vrule#&\strut\;\;#\hfill\cr
height3pt&\omit&&\omit&&\omit&&\omit&\cr
& \lower5pt\hbox{$\!\!(i_1,i_2,i_3)$} &
& \lower5pt\hbox{$N_G(A_1,A_2,A_3)/|A_1^G|$} && \!\!length &&\!\!proof of &\cr
&\omit&&\omit&& \lower-3pt\hbox{\!of $\Delta$} && \lower-3pt\hbox{\!Th.~\thGU} &\cr
\noalign{\hrule}
height3pt&\omit&&\omit&&\omit&&\omit&\cr
&$(2,2,2)$&& $(2q^2\delta_L\pm q-2)\delta_{111}  \quad[\delta_{111}=1]$        &&\omit && \evi&\cr
&$(3,2,2)$&& $2\delta_L\delta_{111} \quad[\delta_{111}=1]$                     &&\omit && \evi&\cr
&$(3,3,2)$&& $q^2(1\mp\delta_{111})+(q-1)\delta_{111}-4q\delta_L\delta_{111}$  &&\omit && \evi&\cr

&$(3,3,3)$&& $q^2(q^2-2) + q(q^2\pm2q-2)\delta_{111}$                          &&\omit &&\evi&\cr
&$(4,2,2)$&& $0 \quad \text{[false]}$                                          &&\omit &&\evi&\cr
&$(4,3,2)$&& 
             $0 \quad \text{[false]}$                                          &&\omit &&\evi&\cr
&$(4,3,3)$&& 
             $q(q\pm1)^2(q\mp1)$                                               &&\omit &&\evi&\cr
&$(4,4,2)$&& $2(q^2-1)\delta_L\delta_{111} \quad[\delta_{111}=1]$              &&\omit &&\evi&\cr
&$(4,4,3)$&& $(q\pm1)^2(q\mp1)\delta_{111} \quad[\delta_{111}=1]$              &&\omit &&\evi&\cr
height3pt&\omit&&\omit&&\omit&&\omit&\cr
&$(4,4,4)$&& $(2q^2\delta_L\pm1)\delta_{111} + q(q\mp1)\delta_{112}\delta_{121}\delta_{211}$
                                                                               &&\omit &&\omit&\cr 
&\omit   && $\qquad\qquad\qquad\qquad\qquad[\delta_{111}+\delta_{112}\delta_{121}\delta_{211}=1]$
                                                                               &&\omit &&\evi&\cr
height3pt&\omit&&\omit&&\omit&&\omit&\cr
&$(5,2,2)$&& $q\delta_{211} \quad[\delta_{211}=1]$                             &&\omit &&\evi&\cr
&$(5,3,2)$&& 
             $q(q\pm1)(1\mp\delta_{211})$                                      &&\omit &&\TT &\cr
&$(5,3,3)$&& $q(q\mp1)^2\big((q\mp2) + \delta_{211}\big)$                      &&\omit &&\evi&\cr
&$(5,4,2)$&& $(q\mp q\pm1)\delta_{111}+q\delta_{211} \quad[\delta_{111}+\delta_{211}=1]$
                                                                               &&\omit &&\evi&\cr
&$(5,4,3)$&& $(q\pm1)\big(q+(2q\delta_L-1)\delta_{111}\mp q\delta_{211}\big)$  &&\omit &&\TT &\cr
&$(5,4,4)$&& $(q\pm1)\delta_{111}+q\delta_{211}(1\mp\delta_{112}\delta_{121})
   \quad[\delta_{111}+\delta_{211}=1]$                                         &&\omit &&\TT&\cr
&$(5,5,2)$&& $q^2\!\pm q + 2\big((q\!-\!1)^2\delta_L\!-1\big)\delta_{111} \mp (q^2\!\pm q)\Delta_1\mp q^2\delta_{221}$
                                                                               &&\omit &&\TT&\cr
height3pt&\omit&&\omit&&\omit&&\omit&\cr
&$(5,5,3)$&& $(q\pm1)\big(q(q^2\mp2q-2)+(q^2-4q\delta_L+1)\delta_{111}$        &&\omit &&\omit &\cr 
&$       $&& $\qquad\qquad\qquad\qquad\qquad\qquad   + q(q\pm1)\Delta_1
                                                     +q^2\delta_{221}\big)$    &&\omit &&\evi &\cr
height3pt&\omit&&\omit&&\omit&&\omit&\cr
&$(5,5,4)$&&
             $q(q\pm1)\big(\delta_{112}\delta_{121}\delta_{211}
               \mp\delta_{121}\mp\delta_{211}+1\big)$                           &&\omit &&\omit&\cr 
&$       $&& $\qquad\qquad\qquad\qquad
                \mp q^2\delta_{221}+(2q^2\delta_L-2q\mp1)\delta_{111}$          &&\omit &&\TT &\cr
height3pt&\omit&&\omit&&\omit&&\omit&\cr
&$(5,5,5)$&& $q(q\pm1)(q^2\mp3q-2+q\Delta_1) + (q^3\pm3q^2-2q^2+3q$             &&\omit &&\omit&\cr 
&$       $&& $\quad
              \pm1)\delta_{111}+q(q\pm1)^2\big(\Delta_2\mp\delta_{112}\delta_{121}\delta_{211}\big)
              +q^3\delta_{222}$ &&\omit && \evi&\cr
height3pt&\omit&&\omit&&\omit&&\omit&\cr
&$(6,2,2)$&& $(q\pm1)\Delta_0
                                                            \quad[\Delta_0=1]$  &&\omit &&\evi&\cr
&$(6,3,2)$&& $(q\pm1)^2(1\mp\Delta_0)$
                                                                                &&\omit &&\TT &\cr
&$(6,3,3)$&& $(q\pm1)^2\big(q^2\mp2q-1+(q\pm1)
                \Delta_0\big)$ &&\omit &&\evi&\cr
&$(6,4,2)$&& $(q\pm1)\Delta_1 \quad[\Delta_1=1]$
     &&\omit &&\evi &\cr
&$(6,4,3)$&& $(q\pm1)^2(1\mp\Delta_1)$&&\omit &&\TT &\cr 
&$(6,4,4)$&& $(q\pm1)\Delta_2
                     \mp q\Delta \quad[\Delta_2=1]$&& \;\;6 &&\evi&\cr
&$(6,5,2)$&& $(q\pm1)\big((q\pm1)(1\mp\Delta_1)\mp q\Delta_0\big)$              &&\omit &&\TT &\cr
&$(6,5,3)$&& $(q\pm1)^2\big( (q^2\mp3q-1) + (q\pm1)\Delta_1 + q\Delta_0\big)$   &&\omit &&\evi &\cr
&$(6,5,4)$&& $(q\pm\!1)\big((q\pm\!1)(1\!\mp\!\Delta_2)
      \mp q(\delta_{121}\!+\!\delta_{221}\!+\!\delta_{321}) + q\Delta\big)$     && \;\;6 &&\TT &\cr
height3pt&\omit&&\omit&&\omit&&\omit&\cr
&$(6,5,5)$&& $(q\pm1)\big( (q\pm1)(q^2\mp4q-1) + (q\pm1)^2\Delta_2$             &&\omit &&\omit &\cr
&$       $&& $\qquad\qquad \qquad\qquad
       +q(q\pm1)(\Delta_1\mp\Delta) + q^2\Delta_0\big)$                         && \;\;6 &&\STP &\cr
height3pt&\omit&&\omit&&\omit&&\omit&\cr
&$(6,6,2)$&& $(q\pm1)\big((q\pm1) \mp q\Delta_0 +(2q\mp1)\Delta\big)$           && \;\;6 &&\TT &\cr
&$(6,6,3)$&& $q(q\pm1)^2\big( q\mp4 +\Delta_0\mp\Delta\big)$                    && \;\;6 &&\STP &\cr
&$(6,6,4)$&& $(q\pm1)\big(1 + q(1\mp\Delta_1)\big) + q^2\Delta$                 && \;18 &&\STP &\cr
&$(6,6,5)$&& $q(q\pm1)\big((q\pm1)(q\mp5)+(q\pm1)\Delta_1 + q\Delta_0\mp q\Delta\big)$
                                                                                && \;18 &&\STP &\cr
&$(6,6,6)$&& $(q\pm1)^2(q^2\mp6q+1)+q^2(q\pm1)\Delta_0\mp q^3\Delta$            && \;36 &&\STP &\cr
height3pt&\omit&&\omit&&\omit&&\omit&\cr
}\hrule}
\endinsert

\midinsert
\noindent Table \tabTripleProd\ (continued-1)
\smallskip
\vbox{\offinterlineskip
\hrule
\halign{&\vrule#&\strut\quad#\hfill\cr
height3pt&\omit&&\omit&&\omit&&\omit&\cr
& \lower5pt\hbox{$\!\!(i_1,i_2,i_3)$} &
& \lower5pt\hbox{$N_G(A_1,A_2,A_3)/|A_1^G|$} && \!\!\!\!\!length && proof of &\cr
&\omit&&\omit&& \lower-3pt\hbox{\!\!of $\Delta$} && \lower-3pt\hbox{Th.~\thGU} &\cr
\noalign{\hrule}
height3pt&\omit&&\omit&&\omit&&\omit&\cr
&$(7,2,2)$&& $(q\mp1)\delta_{111}  \quad[\delta_{111}=1]$&& &&\evi&\cr
&$(7,3,2)$&& $(q^2-1)(1\mp\delta_{111})$&& &&\TT&\cr
&$(7,3,3)$&& $(q\pm1)(q^2-1)\big(q\mp1 + \delta_{111}\big)$&& &&\evi&\cr
&$(7,4,2)$&& $(q\mp1)\delta_{111}\quad [\delta_{111}=1]$ && &&\evi&\cr
&$(7,4,3)$&& $(q^2-1)(1\mp\delta_{111})$&& &&\TT &\cr
&$(7,4,4)$&& $(q\mp1)\delta_{111}\quad [\delta_{111}=1]$&& &&\evi&\cr
&$(7,5,2)$&& $(q\mp1)\big( (q\pm1)(1\mp\delta_{111}) \mp q\delta_{121}\big)$&& &&\TT &\cr
&$(7,5,3)$&& $(q^2-1)\big(q^2\mp q+1 + (q\pm1)\delta_{111} + q\delta_{121}\big)$&& &&\evi&\cr
&$(7,5,4)$&& $(q\mp1)\big((q\pm1)(1\mp\delta_{111})\mp q\delta_{121}\big)$&& &&\TT &\cr
&$(7,5,5)$&& $(q\mp1)\big((q\pm1)(q^2\mp2q-1) + (q\pm1)^2\delta_{111}$ && && &\cr
&$       $&& $\qquad\qquad\qquad\qquad
                 + q(q\pm1)\Delta_1 + q^2\delta_{122}\big)$&& &&\evi&\cr
&$(7,6,2)$&& $(q\mp1)\big(q\pm1 \mp q\Delta_0\big)$&& &&\evi&\cr
&$(7,6,3)$&& $q(q^2-1)\big( q\mp2 + \Delta_0\big)$&& &&\evi&\cr
&$(7,6,4)$&& $(q\mp1)\big(q\pm1 \mp q\Delta_1\big)$&& &&\evi&\cr
&$(7,6,5)$&& $q(q\mp1)\big( (q\pm1)(q\mp3) + (q\pm1)\Delta_1 + q\Delta_0\big)$&& &&\evi&\cr
&$(7,6,6)$&& $(q\mp1)\big((q\pm1)(q^2\mp4q+1) + q^2\Delta_0\big)$&& &&\evi&\cr
&$(7,7,2)$&& $(q\mp1)\big(1 + q(1\mp\delta_{111}) \pm \Delta\big)$  && 2 &&\evi&\cr
&$(7,7,3)$&& $q(q^2-1)\big( q + \delta_{111} \pm
            \Delta\big)$&& 2 &&\evi&\cr
&$(7,7,4)$&& $(q\mp1)\big(q\pm1\mp q\delta_{111}\big) + q^2\Delta$ && 2 &&\evi&\cr
&$(7,7,5)$&& $q(q\mp1)\big( q^2-1 + (q\pm1)\delta_{111} 
           + q\delta_{112} \pm q\Delta\big)$ && 2 &&\evi&\cr
&$(7,7,6)$&& $(q\mp1)\big( (q^2-1)(q\mp1) + q^2\Delta_0\big)$&& && \evi&\cr
&$(7,7,7)$&& $(q^4-1) + q^2(q\mp1)\delta_{111} \pm q^3\Delta$  && 4 &&\evi&\cr
&$(8,2,2)$&& $0 \quad\text{[false]}$ && &&\evi&\cr
&$(8,3,2)$&& $q^2\mp q+1$&& &&\evi&\cr
&$(8,3,3)$&& $(q^2\mp q+1)(q^2\pm q-1)$&& &&\evi&\cr
&$(8,4,2)$&& $0\quad\text{[false]}$&& &&\evi&\cr
&$(8,4,3)$&& $q^2\mp q+1$&& &&\evi&\cr
&$(8,4,4)$&& $0\quad\text{[false]}$&& &&\evi&\cr
&$(8,5,2)$&& $q^2\mp q+1$&& &&\evi&\cr
&$(8,5,3)$&& $(q^2-1)(q^2\mp q+1)$&& &&\evi&\cr
&$(8,5,4)$&& $q^2\mp q+1$&& &&\evi&\cr
&$(8,5,5)$&& $(q^2\mp q+1)(q^2\pm q-1)$&& &&\evi&\cr
&$(8,6,2)$&& $q^2\pm q+1$&& &&\evi&\cr
&$(8,6,3)$&& $q(q\mp1)(q^2\mp q+1)$&& &&\evi&\cr
&$(8,6,4)$&& $q^2\mp q+1$&& &&\evi&\cr
&$(8,6,5)$&& $q(q\mp2)(q^2\mp q+1)$&& &&\evi&\cr
&$(8,6,6)$&& $(q^2\mp q+1)(q^2\mp 3q+1)$&& &&\evi&\cr
&$(8,7,2)$&& $q^2\mp q+1$&& &&\evi&\cr
&$(8,7,3)$&& $q(q^3\pm 1)$&& &&\evi&\cr
&$(8,7,4)$&& $q^2\mp q+1$&& &&\evi&\cr
&$(8,7,5)$&& $q^2(q^2\mp q+1)$&& &&\evi&\cr
&$(8,7,6)$&& $(q^2\mp q+1)^2$&& &&\evi&\cr
&$(8,7,7)$&& $q^4+q^2+1$ && &&\evi&\cr
height3pt&\omit&&\omit&&\omit&&\omit&\cr
}\hrule}
\endinsert

\midinsert
\noindent Table \tabTripleProd\ (continued-2)
\smallskip
\vbox{\offinterlineskip
\hrule
\halign{&\vrule#&\strut\quad#\hfill\cr
height3pt&\omit&&\omit&&\omit&&\omit&\cr
& \lower5pt\hbox{$\!\!(i_1,i_2,i_3)$} &
& \lower5pt\hbox{$N_G(A_1,A_2,A_3)/|A_1^G|$} && \!\!\!\!\!length && proof of &\cr
&\omit&&\omit&& \lower-3pt\hbox{\!\!of $\Delta$} && \lower-3pt\hbox{Th.~\thGU} &\cr
\noalign{\hrule}
height3pt&\omit&&\omit&&\omit&&\omit&\cr
&$(8,8,2)$&& $(q^2\mp q+1)(1-\Delta'/3)$ && $9$ &&\TT &\cr
&$(8,8,3)$&& $q(q^2\mp q+1)\big( q\pm2 \mp \Delta'/3)$ && $9$ &&\evi&\cr
&$(8,8,4)$&& $q^2\mp q+1$ &&\omit &&\evi&\cr
&$(8,8,5)$&& $q(q^3\pm1)$ &&\omit &&\evi&\cr
&$(8,8,6)$&& $(q^2+1)(q^2\mp q+1)$ &&\omit &&\evi&\cr
&$(8,8,7)$&& $(q\pm 1)(q^3\pm 1)$ &&\omit &&\evi&\cr
&$(8,8,8)$&& $(q^2\mp q+1)(q^2\pm 3q+1) \mp q^3\Delta'/3$ && $27$ &&\STP&\cr
height3pt&\omit&&\omit&&\omit&&\omit&\cr
}\hrule}
\endinsert

We do the following substitutions
(we may do them because of the determinant relation):
\roster
\item"$(i)$"
       $\delta_{\vec a}^2=\delta_{\vec a}$;
\item"$(ii)$"
      $\delta_{\vec a}\delta_{\vec b}=0$ if $\vec a$ and $\vec b$
      differ at exactly one position, i.~e., if there exists $\nu_0$
      such that $a_\nu=b_\nu$ if and only if $\nu=\nu_0$,
      for example,
      $\delta_{122}\delta_{132}=0$;
\item"$(iii)$"
      $\delta_{\vec a}=0$ if there exists $\nu_0$ such that
      $a_\nu\le n'_\nu$ if and only if $\nu=\nu_0$, for
      example, we set $\delta_{321}=0$ if $(i_1,i_2,i_3)=(7,5,4)$;
\item"$(iv)$"
      $\delta_{111}\delta_{n_{i_1},n_{i_2},n_{i_3}}
           =\delta_{111}$ if $i_1,i_2,i_3\le 5$;
\item"$(v)$"
      $\delta_{111}=0$ if $i_1\in\{4,5\}$ and $\{i_2,i_3\}\subset\{2,3\}$.
\endroster
The result of computation is presented in Table \tabTripleProd.
Recall that $\delta_L$ is defined by (\eqDeltaL).
In the third column, which is entitled ``length of $\Delta$'', we give the number
of monomials in $\Delta$ or in $\Delta'$ survived after the substitutions
$(i)$--$(v)$.
If there are restrictions on $\delta_{\vec a}$ imposed by the rank
condition, then we write them in the brackets in the second column
(if the rank condition is never satisfied,
then we write ``[false]'').

It is clear from Table \tabTripleProd\
that $N_G(A_1,A_2,A_3)=0$ in the cases $(i)$--$(ix)$ of Theorem \thGU(a).

Also, when $G=GL$, it is clear from Table \tabTripleProd\
that $N_G(A_1,A_2,A_3)\ne0$ unless the cases $(viii)$ and $(ix)$ of Theorem \thGU;
maybe, it worth to
note only that $\Delta\le\delta_{1,1,n_{i_3}}$ for $i_1=i_2=7$, $i_3\in\{2,3,4,5\}$ and that
for $(i_1,i_2,i_3)=(8,8,2)$ the proof is the same as in the case $G=GU$.

In the last column we give a reference to a proof of Theorem \thGU(b) for $G=GU$ and $q\ge5$
in the corresponding case
(``ev.'' means ``evident'').
The case of $G=GU$, $q=2,3,4$, is done in \S\sectSmallq\ and \S\sectQthreeQfour.

Table \tabTPExcept\ serves to prove Theorem \thGU(b) for the triples
$(i_1,i_2,i_3)$ appearing in cases
$(i)$, $(ii)$, $(iii)$, $(v)$ of Theorem \thGU(a).
In the second column we write condition $(*)$ on $\delta_{\vec a}$.
It is a condition which is equivalent to the fact that the hypothesis
of Theorem \thGU(b) is satisfied, i.~e., the conditions $(i)$--$(v)$ are not
satisfied for any permutation of $(i_1,i_2,i_3)$ and for any renumbering of
the eigenvalues under (\eqEV).
As in Table \tabTripleProd,
the rank condition is written in the brackets.
In the third column we write the structure constant for $G=GU$ under
condition $(*)$.
In each case it is obviously nonzero for $q\ge5$.

\midinsert
\noindent Table \tabTPExcept.
\smallskip
\vbox{\offinterlineskip
\hrule
\halign{&\vrule#&\strut\quad#\hfill\cr
height3pt&\omit&\omit&\omit&\omit&\omit&\cr
& $\!\!(i_1,i_2,i_3)$ &\omit
& condition $(*)$ &\omit& $N_G(A_1,A_2,A_3)/|A_1^G|$ under $(*)$ for $G=GU(3,q^2)$ &\cr
height3pt&\omit&\omit&\omit&\omit&\omit&\cr
\noalign{\hrule}
height3pt&\omit&\omit&\omit&\omit&\omit&\cr
& $(5,3,2)$ &\omit& $\delta_{211}=0$ &\omit& $q(q+1)$ &\cr
& $(5,4,3)$ &\omit& $\delta_{211}=0$ &\omit& $(q+1)(q-\delta_{111})$ &\cr
& $(5,4,4)$ &\omit& $\delta_{112}\delta_{121}\delta_{211}=0$ &\omit&
       $(q+1)\delta_{111} + q\delta_{211}\quad [\delta_{111}+\delta_{211}=1]\;\;$ &\cr
& $(5,5,2)$ &\omit& $\Delta_1=0$ &\omit& $q(q+1) - 2\delta_{111} - q^2\delta_{221}$ &\cr
& $(5,5,4)$ &\omit& $\delta_{211}=\delta_{121}=0$ &\omit&
        $q(q+1) - q^2\delta_{221} - (2q+1)\delta_{111}$ &\cr
& $(6,3,2)$ &\omit& $\Delta  =0$ &\omit& $(q+1)^2$ &\cr
& $(6,4,3)$ &\omit& $\Delta_1=0$ &\omit& $(q+1)^2$ &\cr
& $(6,5,2)$ &\omit& $\Delta_1=0$ &\omit& $(q+1)\big( 1 + q(1-\Delta_0)\big)$ &\cr
& $(6,5,4)$ &\omit& $\Delta_2=0$ &\omit&
        $(q+1)\big(1+q(1+\Delta-\delta_{121}-\delta_{221}-\delta_{321})\big)$ &\cr
& $(6,6,2)$ &\omit& $\Delta=0$ &\omit& $(q+1)\big( 1 + q(1-\Delta_0)\big)$ &\cr
& $(7,3,2)$ &\omit& $\delta_{111}=0$ &\omit& $q^2-1$ &\cr
& $(7,4,3)$ &\omit& $\delta_{111}=0$ &\omit& $q^2-1$ &\cr
& $(7,5,2)$ &\omit& $\delta_{111}=0$ &\omit& $(q-1)\big( 1 + q(1-\delta_{121})\big)$ &\cr
& $(7,5,4)$ &\omit& $\delta_{111}=0$ &\omit& $(q-1)\big( 1 + q(1-\delta_{121})\big)$ &\cr
& $(8,8,2)$ &\omit& $\Delta'=0$ &\omit& $q^2-q+1$ &\cr
height3pt&\omit&\omit&\omit&\omit&\omit&\cr
}\hrule}
\endinsert


\subhead\sectQthreeQfour. The cases of $GL(3,2)$ and $GU(3,q^2)$ for $q=3,4$
\endsubhead

These cases are treated in [\refKarni]:
p.~64 for $GL(3,2)$, pp.~69--71 for $GU(3,3^2)$ and pp.~89--93 for $GU(3,4^2)$.
The correspondence between the notation of conjugacy classes
in [\refEnnola] (used in this paper) and the notation in [\refKarni]
is given in Tables \tabSLthreeQtwo,
\tabSUthreeQthree\ and \tabSUthreeQfour. Note that in all these cases $3$ does not
divide $q\pm1$, hence it is enough to consider the case of $SU$ instead of $GU$.

\midinsert
\noindent Table \tabSLthreeQtwo. Notation correspondence for conjugacy classes in $SL(3,2)=GL(3,2)$
\smallskip
\vbox{\offinterlineskip
\hrule
\halign{&\vrule#&\strut\quad#\;\hfill\cr
height3pt&\omit&\omit&\omit&&\omit&\omit&\omit&&\omit&\omit&\omit&\cr
& in [\refKarni] &\omit& in \S\sectCCG\;\; &
& in [\refKarni] &\omit& in \S\sectCCG &
& in [\refKarni] &\omit& in \S\sectCCG &\cr
height3pt&\omit&\omit&\omit&&\omit&\omit&\omit&&\omit&\omit&\omit&\cr
\noalign{\hrule}
height3pt&\omit&\omit&\omit&&\omit&\omit&\omit&&\omit&\omit&\omit&\cr
& 1A &\omit& $C_1^{(0)}$ && 3B &\omit& $C_7^{(1)}=C_7^{(2)}$ && 7A &\omit& $C_8^{(1)}=C_8^{(2)}=C_8^{(4)}$ &\cr
height1pt&\omit&\omit&\omit&&\omit&\omit&\omit&&\omit&\omit&\omit&\cr
& 2A &\omit& $C_2^{(0)}$ && 4B &\omit& $C_3^{(0}$            && 7A &\omit& $C_8^{(3)}=C_8^{(5)}=C_8^{(6)}$ &\cr
height1pt&\omit&\omit&\omit&&\omit&\omit&\omit&&\omit&\omit&\omit&\cr
height1pt&\omit&\omit&\omit&&\omit&\omit&\omit&&\omit&\omit&\omit&\cr
}\hrule}
\endinsert

\midinsert
\noindent Table \tabSUthreeQthree. Notation correspondence for conjugacy classes in $SU(3,3^2)$
\smallskip
\vbox{\offinterlineskip
\hrule
\halign{&\vrule#&\strut\quad#\hfill\cr
height3pt&\omit&\omit&\omit&&\omit&\omit&\omit&&\omit&\omit&\omit&\cr
& in [\refKarni] &\omit& in [\refEnnola] &
& in [\refKarni] &\omit& in [\refEnnola] &
& in [\refKarni] &\omit& in [\refEnnola] &\cr
height3pt&\omit&\omit&\omit&&\omit&\omit&\omit&&\omit&\omit&\omit&\cr
\noalign{\hrule}
height3pt&\omit&\omit&\omit&&\omit&\omit&\omit&&\omit&\omit&\omit&\cr
& 1A &\omit& $C_1^{(0)}$     && 4B &\omit& $C_4^{(3,2)}$   && 8A &\omit& $C_7^{(1,1)}=C_7^{(1,5)}$ &\cr
height1pt&\omit&\omit&\omit&&\omit&\omit&\omit&&\omit&\omit&\omit&\cr
& 2A &\omit& $C_4^{(2,0)}$   && 4C &\omit& $C_6^{(0,1,3)}$ && 8B &\omit& $C_7^{(3,3)}=C_7^{(3,7)}$ &\cr
height1pt&\omit&\omit&\omit&&\omit&\omit&\omit&&\omit&\omit&\omit&\cr
& 3A &\omit& $C_2^{(0)}$     && 6A &\omit& $C_5^{(2,0)}$   && 12A &\omit& $C_5^{(1,2)}$            &\cr
height1pt&\omit&\omit&\omit&&\omit&\omit&\omit&&\omit&\omit&\omit&\cr
& 3B &\omit& $C_3^{(0)}$ && 7A &\omit& $C_8^{(4)}=C_8^{(8)}=C_8^{(16)}$ && 12B &\omit& $C_5^{(3,2)}$ &\cr
height1pt&\omit&\omit&\omit&&\omit&\omit&\omit&&\omit&\omit&\omit&\cr
& 4A &\omit& $C_4^{(1,2)}$ && 7B &\omit& $C_8^{(12)}=C_8^{(20)}=C_8^{(24)}$ &&\omit&\omit&\omit&\cr
height1pt&\omit&\omit&\omit&&\omit&\omit&\omit&&\omit&\omit&\omit&\cr
}\hrule}
\endinsert

\midinsert
\noindent Table \tabSUthreeQfour. Notation correspondence for conjugacy classes in $SU(3,4^2)$
\smallskip
\vbox{\offinterlineskip
\hrule
\halign{&\vrule#&\strut\quad#\hfill\cr
height3pt&\omit&\omit&\omit&&\omit&\omit&\omit&&\omit&\omit&\omit&\cr
& in [\refKarni] &\omit& in [\refEnnola] &
& in [\refKarni] &\omit& in [\refEnnola] &
& in [\refKarni] &\omit& in [\refEnnola] &\cr
height3pt&\omit&\omit&\omit&&\omit&\omit&\omit&&\omit&\omit&\omit&\cr
\noalign{\hrule}
height3pt&\omit&\omit&\omit&&\omit&\omit&\omit&&\omit&\omit&\omit&\cr
& 1A &\omit& $C_1^{(0)}$   && 5E &\omit& $C_6^{(0,1,4)}$ && 13C &\omit& $C_8^{(20)}=C_8^{(50)}=C_8^{(60)}$ &\cr
height1pt&\omit&\omit&\omit&&\omit&\omit&\omit&&\omit&\omit&\omit&\cr
& 2A &\omit& $C_2^{(0)}$   && 5F &\omit& $C_6^{(0,2,3)}$ && 13D &\omit& $C_8^{(35)}=C_8^{(40)}=C_8^{(55)}$ &\cr
height1pt&\omit&\omit&\omit&&\omit&\omit&\omit&&\omit&\omit&\omit&\cr
& 3A &\omit& $C_7^{(0,5)}
            =C_7^{(0,10)}$ && 10A &\omit& $C_5^{(1,3)}$  && 13A &\omit& $C_8^{( 5)}=C_8^{(15)}=C_8^{(45)}$ &\cr
height1pt&\omit&\omit&\omit&&\omit&\omit&\omit&&\omit&\omit&\omit&\cr
& 4A &\omit& $C_3^{(0)}$   && 10B &\omit& $C_5^{(2,1)}$  && 13B &\omit& $C_8^{(10)}=C_8^{(25)}=C_8^{(30)}$ &\cr
height1pt&\omit&\omit&\omit&&\omit&\omit&\omit&&\omit&\omit&\omit&\cr
& 5A &\omit& $C_4^{(1,3)}$ && 10C &\omit& $C_5^{(4,2)}$  && 15A &\omit& $C_7^{(3,8)}=C_7^{(3,13)}$         &\cr
height1pt&\omit&\omit&\omit&&\omit&\omit&\omit&&\omit&\omit&\omit&\cr
& 5B &\omit& $C_4^{(2,1)}$ && 10D &\omit& $C_5^{(3,4)}$  && 15B &\omit& $C_7^{(1,1)}=C_7^{(1,11)}$         &\cr
height1pt&\omit&\omit&\omit&&\omit&\omit&\omit&&\omit&\omit&\omit&\cr
& 5C &\omit& $C_4^{(4,2)}$ &&\omit&\omit& $           $  && 15C &\omit& $C_7^{(2,2)}=C_7^{(2,7)}$          &\cr
height1pt&\omit&\omit&\omit&&\omit&\omit&\omit&&\omit&\omit&\omit&\cr
& 5D &\omit& $C_4^{(3,4)}$ &&\omit&\omit& $           $  && 15D &\omit& $C_7^{(4,4)}=C_7^{(4,14)}$         &\cr
height3pt&\omit&\omit&\omit&&\omit&\omit&\omit&&\omit&\omit&\omit&\cr
}\hrule}
\endinsert

\subhead\sectTPtwo. Proof of Theorem \thGU\ for $m=3$
\endsubhead
Here we complete the proof for triples $(i_1,i_2,i_3)$ not covered by Table \tabTPExcept.
In this section $G=GU$.

\subhead The case $(i_1,i_2,i_3)=(6,5,5)$
\endsubhead
We have 
$$
\split
  \Delta_1-\Delta&=
   \delta_{112}(1-\delta_{211}\delta_{321})
  +\delta_{212}(1-\delta_{311}\delta_{121})
  +\delta_{312}(1-\delta_{111}\delta_{221})
\\&
  +\delta_{121}(1-\delta_{211}\delta_{312})
  +\delta_{221}(1-\delta_{311}\delta_{112})
  +\delta_{321}(1-\delta_{111}\delta_{212}) \ge 0.
\endsplit
$$

\subhead The case $(i_1,i_2,i_3)=(6,6,3)$
\endsubhead
We have
$$
\split
  \Delta_0 - \Delta &=
   \delta_{111}(1-\delta_{221}\delta_{331}-\delta_{231}\delta_{321})
  +\delta_{121}(1-\delta_{211}\delta_{331}-\delta_{231}\delta_{311})
\\&
  +\delta_{131}(1-\delta_{211}\delta_{321}-\delta_{221}\delta_{311})
  + \sum_{\vec a\in\vec n; a_1>1}\delta_{\vec a} \ge 0
\endsplit
$$

\subhead The case  $(i_1,i_2,i_3)=(6,6,4)$
\endsubhead
If $\Delta>0$, 
then there exist permutations of the eigenvalues such that the
product of corresponding diagonal matrices is the identity matrix.
So, we consider only the case when $\Delta=0$.
In this case
$N_G(A_1,A_2,A_3)/|A_1^G|=(q+1)\big( 1 + q(1-\Delta_1)\big)$ which cannot
be zero for any integers $q>1$ and $\Delta_1$.

\subhead The case   $(i_1,i_2,i_3)=(6,6,5)$
\endsubhead
Here we write for shortness $\nu^\alpha$ instead of
$\alpha(\nu)$. We have
$\Delta=\sum_{\alpha\in S_3} \sum_{\beta\in\Cal A_{6,5}}
      \delta_{1,1^\alpha,1^\beta}\delta_{2,2^\alpha,2^\beta}
       \delta_{3,3^\alpha,3^\beta}=\sum_{\alpha\in S_3}E(\alpha)$ where
$$
    E(\alpha)=
     \delta_{1,1^\alpha,1} \delta_{2,2^\alpha,1} \delta_{3,3^\alpha,2}
   + \delta_{1,1^\alpha,1} \delta_{2,2^\alpha,2} \delta_{3,3^\alpha,1}
   + \delta_{1,1^\alpha,2} \delta_{2,2^\alpha,1} \delta_{3,3^\alpha,1}.
$$
Summating $E(\alpha)$ separately over odd and even permutations $\alpha$
and estimating each triple product of the deltas by one of its factors,
we obtain
$$
\split
  &\sum_{\text{odd $\alpha$}} E(\alpha)\le
   \sum_{\text{odd $\alpha$}} \big(\delta_{3,3^\alpha,2}+\delta_{2,2^\alpha,2}
         + \delta_{1,1^\alpha,2}\big) = \Delta_0,
\\
  &\sum_{\text{even $\alpha$}} E(\alpha)\le
   \sum_{\text{even $\alpha$}} \big(\delta_{1,1^\alpha,1}+\delta_{3,3^\alpha,1}
         + \delta_{2,2^\alpha,1}\big) = \Delta_1
\endsplit
$$
which implies $\Delta_1 + \Delta_0 - \Delta\ge 0$ and the result follows for $q>5$.

Let $q=5$. The above considerations show that the structure constant is positive
when $\Delta_1>0$. So, we suppose that $\Delta_1=0$. Then $\Delta=0$ because each triple
product in $\Delta$ includes some $\delta_{\vec a}$ involved in $\Delta_1$.
If we have two triples of distinct residues mod 6 (the
parameters $(k,l,m)$ of $C_6^{(k,l,m)}$) not of
the same parity, then their pairwise sums attain
all values mod 6 except, maybe one, thus $\Delta_0$ or $\Delta_1$ is nonzero.
So, it remains to consider the case $A_1,A_2\in C_6^{(0,2,4)}$.
In this case, (\eqDet) implies $A_3\in C_5^{(k,l)}$ with $l$ even,
hence $\Delta_0>0$ and the result follows.

\subhead The case  $(i_1,i_2,i_3)=(6,6,6)$
\endsubhead
If $\Delta>0$, 
then there exist permutations of the eigenvalues such that the
product of corresponding diagonal matrices is the identity matrix.
So, we consider only the case when $\Delta=0$.
In this case, the structure constant is positive for
$q>5$ and it is equal to $150\Delta_0-144\ne 0$ for $q=5$.

\subhead The case  $(i_1,i_2,i_3)=(8,8,8)$
\endsubhead

Let the eigenvalues of $A_\nu$ be
$(\lambda_\nu,\lambda_\nu^{q^2},\lambda_\nu^{q^4})$, $\nu=1,2,3$.
Then we have
$$
   \Delta' = \sum_{0\le a,b,c\le 2}\delta_{a,b,c}, \qquad
   \delta_{a,b,c}=\cases 1, & 
   \lambda_1^{q^{2a}}\lambda_2^{q^{2b}}\lambda_3^{q^{2c}}=1,\\
        0, &\text{otherwise}\endcases
$$
It is clear that $\delta_{a,b,c} = \delta_{a',b',c'}$ if
$a-a'\equiv b-b'\equiv c-c'\mod 3$.

We are going to show that there is at most 9 triples $(a,b,c)$ such that 
$\delta_{a,b,c}=1$.
Suppose that one of $\delta_{a,b,c}$ is nonzero.
Without loss of generality we may assume that it is $\delta_{000}$
(otherwise we permute cyclically the eigenvalues of each matrix).
So, we have $\lambda_1\lambda_2\lambda_3=1$.

Let us show that if $\delta_{a,b,c}=1$, then either $a=b=c$
or $a,b,c$ are pairwise distinct (there are only nine such triples).
Suppose that this is not so, say, $a\ne b=c$.             
Then $\delta_{001}\delta_{112}\delta_{220}=1$ or
$\delta_{001}\delta_{112}\delta_{220}=1$ (we consider only the first case).
This means that $\lambda_1\lambda_2\lambda_3^{q^2}=1$.
Combined with $\lambda_1\lambda_2\lambda_3=1$
this yields $\lambda_3^{q^2}=1$, i.e $\lambda_3\in\GF{q^2}$. Contradiction.

Thus, we proved that $\Delta'\le 9$, hence
$$
   N_G(A_1,A_2,A_3)/|A_1^G| 
  \ge (q^2-q+1)(q^2+3q+1) - 3q^3 
%
 = q^4 - q^3 - q^2 + 2q + 1 > 0.
$$


\subhead\sectQP. End of proof of Theorem \thGU\ (the case $m\ge4$)
\endsubhead
Let us prove Theorem \thGU\ for $m=4$.
So, let $m=4$ and let $A_1,\dots,A_4$ be as in Theorem \thGU.

If $G=GL$ and $q\ge3$, then for any $d,\lambda_1,\lambda_2\in\Omega$ there exists
$B\in C_3\cup C_5\cup C_6$ such that $d=\det B$ and $\lambda_1,\lambda_2$ are
eigenvalues. Hence, we can choose $B$ in $C_3\cup C_5\cup C_6$
such that the rank condition is satisfied for both triples $(A_1,A_2,B)$ and
$(B^{-1},A_3,A_4)$. As we have already shown, there are no other restrictions
for triple products in $GL$. This completes the proof of Theorem \thGU\ for $G=GL$.

\proclaim{ Lemma \lemQone }
Let $G=GU$ and $q\ge4$.
Then for any $d,\mu\in\UC$ there exists
$B\in C_7$ such that $\det B=d$ and 
$\lambda_1(B)=\mu$.
\endproclaim

\demo{ Proof } Obvious. \qed\enddemo

\proclaim{ Lemma \lemQtwo }
Let $G=GU$ and $q\ge 5$. Suppose that one of the following conditions holds
\roster
\item"$(i)$"    $\{i_1,i_3\}\not\subset\{2,4\}$;
\item"$(ii)$"   $\{i_1,i_2,i_3\}\subset\{2,4\}$ and $i_4\in\{6,7,8\}$;
\item"$(iii)$"  $i_1=4$, $\{i_2,i_3\}\subset\{2,4\}$, $i_4\in\{3,5\}$;
\item"$(iv)$"   $i_1=i_2=i_3=2$, $i_4\in\{3,5\}$, and $\delta_{1111}=0$;
\item"$(v)$"    $\{i_1,i_2,i_3,i_4\}\subset\{2,4\}$ and $\delta_{1111}=1$;
\item"$(vi)$"   $i_1=i_3=2$, $\{i_2,i_4\}\subset\{2,4\}$ and $\delta_{1111}=0$;
\item"$(vii)$"  $i_1=i_2=i_3=i_4=4$ and $\delta_{1111}=0$.
\endroster
Then $I\in A_1^G\dots A_4^G$.
\endproclaim

\demo{ Proof } We set $d=\det(A_1 A_2)=\det(A_3^{-1} A_4^{-1})$,
$\mu_1=\lambda_1(A_1)\lambda_1(A_2)$, and
$\mu_2=\lambda_1(A_3^{-1})\lambda_1(A_4^{-1})$. We consider the cases $(i)$--$(vii)$
one by one and in each case we find $B$ such that $B\in A_1^G A_2^G$ and
$B^{-1}\in A_3^G A_4^G$. When we choose $B$ in $C_7$, we use Lemma \lemQone.
\medskip

$(i)$.
We choose $B\in C_7$ such that
$\det B=d$ and $\lambda_1(B)\not\in\{\mu_1,\mu_2\}$.

\medskip
$(ii)$.
We choose $B\in C_7$ such that
$\det B=d$ and $\lambda_1(B)=\mu_1$.

\medskip
$(iii)$.
We consider two cases.

\smallskip
{\it Case 1.} $\delta_{1111}=1$, i.~e., $\mu_1=\mu_2$.
We choose $B\in C_3\cup C_5$ such that $\det B=d$ and
$\lambda_1(B)=\mu_1=\mu_2$.

\smallskip
{\it Case 2.} $\delta_{1111}=0$, i.~e., $\mu_1\ne\mu_2$.
Then we choose $B\in C_7$ such that $\det B=d$
and $\lambda_1(B)=\mu_1$.

\medskip
$(iv)$. The choice of $B$ is the same as for $(iii)$, Case 2.

\medskip
$(v)$.
Since $\delta_{1111}=1$, we have $\mu_1=\mu_2$.
So, we choose $B\in C_7$ such that $\det B=d$ and $\lambda_1(B)=\mu_1=\mu_2$.

\smallskip
$(vi)$.
Since $\delta_{1111}=0$, we have $\mu_1\ne\mu_2$.
We choose $B\in C_5\cup C_6$ such that $\det B=d$ and
$\mu_1,\mu_2$ are eigenvalues of $B$.

\smallskip
$(vii)$.
Since $\delta_{1111}=0$, we have $\mu_1\ne\mu_2$.
We choose $B\in C_4\cup C_6$ such that $\det B=d$ and
$\mu_1,\mu_2$ are eigenvalues of $B$.
\qed
\enddemo

For the cases not covered by Lemma \lemQtwo\ we compute the structure
constant in $G=GU$:
$$
\xalignat 3
&(i_1,i_2,i_3,i_4) && \delta_{1111} && N_G(A_1,A_2,A_3,A_4)/|A_1^G| \\
&(3,2,2,2)         &&  1            &&  0 \\
&(5,2,2,2)         &&  1            && (q+3)(q^2-1) \\
&(4,4,4,2)         &&  0            && q(q^2-1)\big(q+1 
                                      -q(\delta_{1121}+\delta_{1211}+\delta_{2111})\\
& && && \qquad\qquad\qquad        +(2q-1)\delta_{1121}\delta_{1211}\delta_{2111}\big)
\endxalignat
$$
This completes the proof of Theorem \thGU\ for $m=4$.

\smallskip
Let $m=5$, $q\ge 5$.
Easy to see that there exists
$B\in(A_1^G A_2^G)\cap(C_3\cup C_5\cup C_6\cup C_7\cup C_8)$.
Then $I\in B^G A_3^G A_4^G A_5^G$. Theorem \thGU\ is proven.


\head\sectS.  Products of conjugacy classes in $SU(3,q^2)$ and $SL(3,q)$. Proof of Theorem \thSU
\endhead

\subhead\sectCTS. The character table of $SU(3,q^2)$ and $SL(3,q)$
\endsubhead
Let $G$ be $GU(3,q^2)$ or $GL(3,q)$ and let $S=\{A\in G\,|\,\det A=1\}$. So,
$S$ is $SU(3,q^2)$ or $SL(3,q)$. The character table of $S$ is computed in [\refSF].
It has some mistakes which are corrected in [\refGAP] (it is
written in the comments in [\refGAP] that the character table for
$SU(3,q^2)$ is taken from [\refGeck]).
Since $G=S\times\UC$ when $3$ does not divide $q\pm1$, we consider
only the case when $q=3r\mp1$.

The conjugacy classes of $S$ are as follows.
Each of $C_3^{(k)}$, $k=0,r,2r$, splits into three classes
$C_3^{(k,l)}$, $l=0,1,2$. The class $C_3^{(k,l)}$ in $SU(3,q^2)$ (resp. in $SL(3,q)$)
consists of matrices which are conjugate in $SL(3,q^2)$ (resp. in $SL(3,q)$)
to\footnote{there is a misprint here in [\refSF].}
$$
   \left(\matrix \omega^k & 0        & 0       \\
                 z^l   & \omega^k & 0       \\
                 0        & 1        & \omega^k
   \endmatrix\right),
     \qquad  z = \cases \rho, & S=SU(3,q^2),\\ \omega, & S=SL(3,q).\endcases
$$
Other conjugacy classes of $G$ contained in $S$ are conjugacy
classes of $S$.

The irreducible characters of $S$ can be described as follows.
We consider the action of the cyclic group of order $q\pm1$ on $\Irr(G)$
such that the action of the generator is
$$
   \chi_{d_j}^{(t)}\mapsto \chi_{d_j}^{(t+1)} \;(j=1,2,3);\quad
   \chi_{d_j}^{(t,u)}\mapsto \chi_{d_j}^{(t+1,u+1)}\;(j=4,5);
$$
$$
   \chi_{d_6}^{(t,u,v)}\mapsto\chi_{d_6}^{(t+1,u+1,v+1)};\quad
   \chi_{d_7}^{(t,u)}\mapsto\chi_{d_7}^{(t+1,u\mp q+1)};\quad
   \chi_{d_8}^{(t)}\mapsto\chi_{d_8}^{(t+q^2\mp q+1)}.
$$
Then the restriction of all characters to $S$ are constant on each 
orbit of this action. All orbits but three are of length $q\pm1$
and their representatives restricted to $S$ are irreducible.
There are three orbits of length $r$, namely the orbits of $\chi_{d_6}^{(0,r,2r)}$
and $\chi_{d_8}^{(u(q^2\mp q+1)/3)}$, $u=1,2$.
Being restricted to $S$, each of these three characters splits into three
irreducible characters. This yields irreducible characters
$\chi_{d_6/3}^{(t)}$, $\chi_{d_8/3}^{(t,u)}$, $t=0,1,2$, $u=1,2$,
such that $\chi_{d_6/3}^{(t)}(A) = {1\over3}\chi_{d_6}^{(0,r,2r)}(A)$ and
$\chi_{d_8/3}^{(t,u)}(A) = {1\over3}\chi_{d_8}^{(u(q^2\mp q+1)/3)}(A)$ when
$A\not\in C_3$. For $A\in C_3^{(k,l)}$, $k,lr\in\{0,r,2r\}$, we have 
$$
  \chi_{d_6/3}^{(t)}(A) = \cases
      q-r,  & l=t,\\
      \quad -r,  & l\ne t,
  \endcases
\qquad
   \chi_{d_8/3}^{(t,u)}(A) = \eps^{uk} \chi_{d_6/3}^{(t)}(A).
$$
where $\eps=f(\omega)=\exp(2\pi i/(q\pm1))$.

Thus, for any function $E$ on $\Irr(S)$, we have
$$
\split
   \sum_{\chi\in\Irr(S)}\!\!\!\!\! E(\chi)
  = {1\over q\pm1}\Big( \sum_{\chi\in\Irr(G)}\!\!\!\!\! E(\chi|_S)\Big)
     &- {1\over 3}\Big( E\big(\chi_{d_6}^{(0,r,2r)}|_S\big)
              + \sum_{u=1}^2 E\big(\chi_{d_8}^{u(q^2\mp q+1)/3}|_S\big)\Big)
\\
     &+ \sum_{t=0}^2\Big( E\big(\chi_{d_6/3}^{(t)}\big) +
           \sum_{u=1}^2 E\big(\chi_{d_8/3}^{(t,u)}\big)\Big)
\endsplit
$$


\subhead\sectSCS. Structure constants for $SU(3,q^2)$ and $SL(3,q)$
\endsubhead
Let $A_1,\dots, A_m\in S$, $A_\nu\in C_{i_\nu}$, $\nu=1,\dots,m$. We suppose that
$i_1=\dots=i_n=3$ and $i_\nu\ne3$ for $\nu>n$.
Let $A_\nu\in C_3^{(k_\nu,l_\nu)}$ for $\nu=1,\dots,n$.

We denote $E_1(\chi)=\chi(A_1)\dots\chi(A_n)$,
$E_2(\chi)=\chi(A_{n+1})\dots\chi(A_m)$, and $E(\chi)=E_1(\chi)E_2(\chi)/\chi(I)^{m-2}$.
Combining the formulas from the previous section with the fact that
$\chi_{d_6}^{(0,r,2r)}(A_\nu)=\mp1$ and $\chi_{d_8}^{u(q^2\mp q+1)/3)}(A_\nu)=\mp\eps^{k_\nu}$
for $\nu\le n$,
we obtain
$$
  E_1\big(\chi_{d_6}^{(0,r,2r)}\big)=(\mp1)^n, \quad
  E_1\big(\chi_{d_8}^{(u(q^2\mp q+1)/3)}\big)=(\mp1)^n\eps^{(k_1+\dots+k_n)u},
$$
$$
  E_2\big(\chi_{d_6/3}^{(t)}\big)=3^{n-m}E_2\big(\chi_{d_6}^{(0,r,2r)}\big),
  \quad
  E_2\big(\chi_{d_8/3}^{(t,u)}\big)=3^{n-m}E_2\big(\chi_{d_8}^{(u(q^2\mp q+1)/3)}\big),
$$
$$
   E_1\big(\chi_{d_8/3}^{(t,u)}\big)=\eps^{(k_1+\dots+k_n)u}E_1\big(\chi_{d_6/3}^{(t)}\big),
   \quad  \chi_{d_6/3}^{(t)}(I)=d_6/3, \quad \chi_{d_8/3}^{(t,u)}(I)=d_8/3,
$$
and finally,
$$
\split
   \bar N_S(A_1,\dots,A_m) = &{\bar N_G(A_1,\dots,A_m)\over q\pm 1}
    +\Bigg(-{(\mp1)^n\over3} + 3^{n-2}\sum_{t=0}^2 E_1\big(\chi_{d_6/3}^{(t)}\big)\Bigg)
\\
   &\times
          \Bigg({E_2\big(\chi_{d_6}^{(0,r,2r)}\big)\over d_6^{m-2}}
        + 
          \sum_{u=1}^2
             {  \eps^{(k_1+\dots+k_n)u}   E_2\big(\chi_{d_8}^{(u(q^2\mp q+1)/3)}\big)
             \over d_8^{m-2}}\Bigg)
\endsplit
$$
In particular, we see from this formula that if $n=0$ or $n=1$, then $\bar N_G=(q\pm1)\bar N_S$, i.~e.,
we have $(I\in A_1^G\dots A_m^G)\Leftrightarrow(I\in A_1^S\dots A_n^S)$. Indeed, if $n=0$, then
the factor
$\big(-{(\mp1)^n\over3}+\dots\;\big)$ is equal to $-1/3+1/9\,(1+1+1)=0$, and if $n=1$, then it
is equal to $\pm1/3 + 1/3\,\big((q-r)-r-r\big)=0$.
This equivalence also follows immediately
from the fact that $C_3^{(k)}$ are the only classes that split in $S$.

\subhead\sectProofS. Triple products  in $SU(3,q^2)$ and $SL(3,q)$. Proof of Theorem \thSU
\endsubhead
Let $m=3$. It is enough to consider the cases $n=2$ and $n=3$.
We use the following notation in Table \tabSU.
If $n=2$, then we set
$$
     \delta^*=\delta^*(A_1,A_2)=\cases 1, & l_1=l_2,\\0,&l_1\ne l_2.\endcases.
$$
If 
$A_3\in C_8^{((q\pm1)k')}$ (the last line of the table), then we set
$$
     \delta_{111}^* = \cases 1, & k_1+k_2+k'\equiv0\mod q\pm1,\\
          0, &\text{otherwise.} \endcases
$$

\midinsert
\noindent Table \tabSU. Structure constants: $S=SU(3,q^2)$ or $SL(3,q)$, $q=3r\mp1$,
$A_\nu\in C_{i_\nu}$
\smallskip
\vbox{\offinterlineskip
\hrule
\halign{&\vrule#&\strut\quad#\hfill\cr
height3pt&\omit&\omit&\omit&&\omit&\cr
& $\!\!(i_1,i_2,i_3)$ &\omit&\omit&& $N_S(A_1,A_2,A_3)/|A_1^S|$ &\cr
height3pt&\omit&\omit&\omit&&\omit&\cr
\noalign{\hrule}
height3pt&\omit&\omit&\omit&&\omit&\cr
& $(3,3,3)$ &\omit& distinct $l_1,l_2,l_3$ \quad && $qr\big(qr+(2qr\mp q+r)\delta_{111}\big)$ &\cr
height3pt&\omit&\omit&\omit&&\omit&\cr
& $(3,3,3)$ &\omit& $l_1=l_2\ne l_3$ && $qr\big(q(r\mp1)-(qr\mp q-r+1)\delta_{111}\big)$  &\cr
height3pt&\omit&\omit&\omit&&\omit&\cr
& $(3,3,3)$ &\omit& $l_1=l_2=l_3$ && $q\big(q(r^2-1)+(2q(r\mp1)^2+r^2\mp1)\delta_{111}\big)\quad$ &\cr
height3pt&\omit&\omit&\omit&&\omit&\cr
& $(3,3,2)$ &\omit& \omit && $\big(q^2 - (q^2\mp q+1)\delta_{111}\big)\delta^*+
               2qr\delta_L\delta_{111}$ &\cr
height3pt&\omit&\omit&\omit&&\omit&\cr
& $(3,3,4)$ &\omit& \omit && $q^2 \delta^*$ &\cr
height3pt&\omit&\omit&\omit&&\omit&\cr
& $(3,3,5)$ &\omit& \omit && $q^2r(q\mp1\mp3\delta^*+\delta_{211})$ &\cr
height3pt&\omit&\omit&\omit&&\omit&\cr
& $(3,3,6)$ &\omit& $\lambda_1(A_3)^r = \lambda_2(A_3)^r$ && $q^2\big((q-1)r\mp2q\delta^*+r\Delta_0\big)$ &\cr
height3pt&\omit&\omit&\omit&&\omit&\cr
& $(3,3,6)$ &\omit& $\lambda_1(A_3)^r\ne\lambda_2(A_3)^r$ && $q^2\big((q-1)r\mp q(1-\delta^*) + r\Delta_0\big)$ &\cr
height3pt&\omit&\omit&\omit&&\omit&\cr
& $(3,3,7)$ &\omit& \omit && $q^2r(q\mp1+\delta_{111})$ &\cr
height3pt&\omit&\omit&\omit&&\omit&\cr
& $(3,3,8)$ &\omit& \omit && $q^2\big( (q-1)r\pm q(\delta^*+\delta_{111}^*-3\delta^*\delta_{111}^*)\big)$ &\cr
height3pt&\omit&\omit&\omit&&\omit&\cr
}\hrule}
\endinsert

It is clear that if $r>1$, 
then
the structure constants are positive except the case when
$i_3\in\{2,4\}$ and $\delta^*=0$ (note that the case $i_3=6$, $q=5$,
$\lambda_1(A_3)^r=\lambda_2(A_3)^r$ is impossible).
This completes the proof of Theorem \thSU\ for $m=3$.

For $m=4$, the proof is the same as in \sectQP. Moreover, since at least two of $A_1,\dots,A_4$
belong to $C_3$, then only Case $(i)$ of Lemma \lemQtwo\ is to be considered.


\head\sectSmallq. The case $q=2$
\endhead

\subhead\sectSmallqG.
Class products in $GU(3,q^2)$ for $q=2$
\endsubhead
Let $G=GU(3,2^2)$, $S=SU(3,2^2)$. Then $|G|=648$, $|S|=216$.
We have the following conjugacy classes in $G$:
$$
\xalignat2
   &\det(A)=1:      && C_1^{(k)}, C_2^{(k)}, C_3^{(k)} (k=0,1,2),\;C_6^{(0,1,2)},\\
   &\det(A)=\rho:   && C_4^{(0,1)}, C_4^{(2,0)}, C_4^{(1,2)},
                      C_5^{(0,1)}, C_5^{(2,0)}, C_5^{(1,2)}, C_8^{(1)},\\
   &\det(A)=\rho^2: && C_4^{(0,2)}, C_4^{(1,0)}, C_4^{(2,1)},
                      C_5^{(0,2)}, C_5^{(1,0)}, C_5^{(2,1)}, C_8^{(2)}
\endxalignat
$$
We see from Table \tabTripleProd\ that $C_6\cdot C_6=C_1\cup C_6$, hence $H=C_6\cup C_1$
is a normal subgroup of $G$ of order 27.
We have $|G/H|=24$ and $S/H=8$.
The sizes of classes and the orders of their representatives in $G/H$ are:
$$
\xalignat8
&\text{Class:}&C_1&^{(k)}&C_2&^{(k)}&C_3&^{(k)}&C_4&^{(k,l)}&C_5&^{(k,l)}&C_6&^{(0,1,2)}&C_8&^{(k)}\\
&\text{Size:} &  1&    &  9&      &  54&     &  12&       &  36&       &  24&         &  72&\\
&\text{Order
      in $G/H$:}&  1&    &  2&      &   4&     &   3&       &   6&       &   1&         &  3&
\endxalignat
$$
The elements of $C_2$ (resp. $C_3$) represent elements of order 2 (resp. 4) in $S/H$.
Since $|H|=|C_2|=27$ and $|C_3|=162$, it follows that $S/H$ has
one element of order 2 and six elements of order $4$. Therefore, $S/H$ is isomorphic to the
unit quaternionic group $Q=\{\pm1,\pm i,\pm j,\pm k\}$. Since the exact sequence
$1\to S/H\to G/H\overset{\det}\to\longrightarrow\{1,\rho,\rho^2\}$ splits, it follows that 
$G/H$ is isomorphic to a semi-direct
product of $Q$ and $\Bbb Z_3$. We denote it by $F$.
Since $G/H$ has no element of order $12$, this product is not direct,
hence $F$ can be identified with
the group whose elements are $\pm a^m$, $\pm ia^m$, $\pm ja^m$, $\pm ka^m$, $m=0,1,2$, subject
to relations $ia=aj$, $ja=ak$, $ka=ai$, $a^3=1$. We denote $a^2$ by $b$.
The conjugacy classes in $F$ are:
$\{1\}$,
$\{-1\}$,
$i^F=\{\pm i,\pm j,\pm k\}$, 
$a^F=\{a,ia,ja,ka\}$,
$-a^F=\{-a,-ia,-ja,-ka\}$,
$b^F=\{b,-ib,-jb,-kb\}$,
$-b^F=\{-b,ib,jb,kb\}$.
Their pairwise products are:
$$
\matrix
 \{1\} &\{-1\}  & i^F &&  a^F           & -a^F          &&  b^F           & -b^F           \\
\{-1\} & \{1\}  & i^F && -a^F           &  a^F          && -b^F           &  b^F           \\
   i^F &  i^F   &  Q  &&  Qa            &  Qa           &&  Qb            &  Qb            \\
\\
   a^F & -a^F   & Qa  &&  Qb            &  Qb           && \{ 1\}\cup i^F & \{-1\}\cup i^F \\
  -a^F &  a^F   & Qa  &&  Qb            &  Qb           && \{-1\}\cup i^F & \{ 1\}\cup i^F \\
\\
   b^F &-b^F    & Qb  && \{ 1\}\cup i^F &\{-1\}\cup i^F &&  Qa            &  Qa            \\
  -b^F & b^F    & Qb  && \{-1\}\cup i^F &\{ 1\}\cup i^F &&  Qa            &  Qa
\endmatrix
$$
Comparing the class sizes and the orders of their representatives, we easily see that the
correspondence between the classes under the projection $G\to F$ is
$$
\xalignat3
    C_1\cup C_6&\to\{1\}  & C_{48}^{(1)}&\to  a^F & C_{48}^{(2)} &\to  b^F\\
            C_2&\to \{-1\}&    C_5^{(1)}&\to -a^F &    C_5^{(2)} &\to -b^F\\
            C_3&\to i^F
\endxalignat
$$
where $C_{48}^{(k)}=(C_4\cup C_8)\cap G^{(k)}$, $C_5^{(k)}=C_5\cap G^{(k)}$,
and $G^{(k)}=\{A\in G\,|\,\det A=\rho^k\}$,
$k=1,2$.
Thus, the multiplication table for the preimages in $G$ of the conjugacy classes of $F$ is
$$
\matrix
  H & C_2 & C_3 && C_{48}^{(1)} &    C_5^{(1)} && C_{48}^{(2)} &    C_5^{(2)}\\
 C_2& H   & C_3 &&    C_5^{(1)} & C_{48}^{(1)} &&    C_5^{(2)} & C_{48}^{(2)}\\
 C_3& C_3 & S   &&      G^{(1)} &      G^{(1)} &&      G^{(2)} &      G^{(2)}\\
\\
 C_{48}^{(1)} &    C_5^{(1)} & G^{(1)} && G^{(2)} & G^{(2)} &&   H\cup C_3 & C_2\cup C_3\\
    C_5^{(1)} & C_{48}^{(1)} & G^{(1)} && G^{(2)} & G^{(2)} && C_2\cup C_3 &   H\cup C_3\\
\\
 C_{48}^{(2)} &    C_5^{(2)} & G^{(2)} &&   H\cup C_3 & C_2\cup C_3 && G^{(1)} & G^{(1)}\\
    C_5^{(2)} & C_{48}^{(2)} & G^{(2)} && C_2\cup C_3 &   H\cup C_3 && G^{(1)} & G^{(1)}
\endmatrix
$$
The above discussion can be summarized as follows

\proclaim{ Proposition \propSmallqOne }
Let $c=(c_1,\dots,c_m)$ is an unordered $m$-tuple of non-trivial conjugacy classes in $F$
such that $\deg_a c_1 +\dots+\deg_a c_m = 0$. We suppose that $c_1=\dots=c_{2n}=\{-1\}$ and
$(c_{2n+1},\dots,c_m)$ contains at most one occurrence of $\{-1\}$.
Then $1\not\in c_1\dots c_m$ if and only if $(c_{2n+1},\dots,c_m)$
is one of $(\{-1\})$, $(i^F)$,
$(\{-1\},i^F)$, 
$(a^F, -b^F)$, $(-a^F,b^F)$, $(\{-1\},a^F,b^F)$, $(\{-1\},-a^F,-b^F)$.
\endproclaim

\demo{ Proof }
It is enough to check that the product of any three non-trivial conjugacy classes
different from $\{-1\}$ is a coset of $Q$ in $F$.
\qed\enddemo

\proclaim{ Proposition \propSmallqTwo }
Let $A_1,\dots,A_m\in G\setminus C_1$ be such that $\det(A_1\dots A_n)=1$.
Let $A_1\in C_{i_1},\dots,A_m\in C_{i_m}$. Suppose that after removing
any number of $6$'s and an even number of $2$'s from $(i_1,\dots,i_m)$,
we obtain one of
$(2)$, $(3)$, $(2,3)$, $(5,4)$, $(8,5)$, $(4,4,2)$, $(5,5,2)$, $(8,4,2)$, $(8,8,2)$.
Then $I\not\in A_1^G\dots A_n^G$.
\endproclaim

\proclaim{ Proposition \propSmallqThree }
Let $A_1,\dots,A_m\in G\setminus C_1$, $m\ge3$,  be such that $\det(A_1\dots A_m)=1$.
Let $A_1\in C_{i_1},\dots,A_m\in C_{i_m}$.
Suppose that the conditions of Proposition \propSmallqTwo\ are not satisfied.
Suppose also that the rank condition {\rm(\eqRk)} holds and
the conditions $(i)$--$(vii)$ of Theorem \thGU(a) are not satisfied for any permutation of $A_1,\dots,A_m$ and for any renumbering of the eigenvalues under restrictions {\rm(\eqEV)}.

Then $I\not\in A_1^G\dots A_m^G$ if and only if one of the following
cases occurs up to changing the order of $A_j$'s,
multiplication them by scalar or simultaneous
replacing of $A_1,\dots,A_m$ by $A_1^{-1},\dots,A_m^{-1}$.
\roster
\item"$(i)$"
             $m=4$, $A_1,A_2,A_3\in C_4^{(0,1)}$ and $A_4\in C_3^{(1)}$;
\item"$(ii)$"
             $m=4$, $A_1,A_2\in C_4^{(0,1)}$, $A_3\in C_4^{(0,2)}$, and $A_4\in C_5^{(1,0)}$.
\endroster
\endproclaim

\demo{ Proof }
Using the structure constants, we computed the products of all $m$-tuples
of conjugacy classes for $m\le 5$.
So we check that the statement is true for $m\le 5$.
The general case easily follows from the following facts.
\roster
\item
"$\bullet$" $C_6^{(0,1,2)}C_6^{(0,1,2)}=H$;
\item
"$\bullet$" $C_2^{(k_1)}C_2^{(k_2)}=C_1^{(k_1+k_2)}\cup C_6^{(0,1,2)}$ for any $k_1$, $k_2$;
\item
"$\bullet$" $C_2^{(k)}C_6^{(0,1,2)}=C_2$ for any $k$;
\item
"$\bullet$" Let $m=4$ or $5$. 
    If $(i_1,\dots,i_m)$ is not as in Proposition \propSmallqTwo\ and
    $\{i_1,\dots,i_m\}\not\subset\{2,6\}$, then $A_1^G\dots A_m^G$ is a coset of $S$ in $G$
    for any $A_1\in C_{i_1},\dots,A_m\in C_{i_m}$;
\endroster
\qed\enddemo

\subhead\sectSmallqS. Class products in $SU(3,q^2)$ for $q=2$
\endsubhead
There are 16 conjugacy classes in $S$. These are:
$$
   C_1^{(k)}, C_2^{(k)}, C_3^{(k,l)}, C_6^{(0,1,2)}, \qquad k,l=0,1,2.
$$
We have $S/H=Q$ and $S/(H\cup C_2)=Q/\{\pm1\}=\Bbb Z_2\oplus\Bbb Z_2$.
The cosets of $H\cup C_2$ in $S$ are: $H\cup C_2$, $C_3^{(*,0)}$, $C_3^{(*,1)}$, $C_3^{(*,2)}$
where $C_3^{(*,l)}$ stands for $C_3^{(0,l)}\cup C_3^{(1,l)}\cup C_3^{(2,l)}$.

\proclaim{ Proposition \propSmallqSU }
Let $A_1,\dots,A_m\in S\setminus C_1$, $m\ge3$, $A_\nu\in C_{i_\nu}$, $\nu=1,\dots,m$.
If $3\in\{i_1\dots,i_m\}$, then $A_1^S\dots A_m^S$ is a coset of $H\cup C_2$ in $S$.
Otherwise $A_1^S\dots A_m^S$ is a coset of $H$ in $H\cup C_2$.
\endproclaim

\demo{ Proof } It is enough to compute the structure constants for all triples $A_1,A_2,A_3\in S$.
\enddemo

\proclaim{ Corollary \corSmallqSU }
Let $A_1,\dots,A_m\in S\setminus C_1$, $m\ge3$, $A_\nu\in C_{i_\nu}$, $\nu=1,\dots,m$.
Then $I\in A_1\dots A_m$ if and only if none of the following conditions holds:
\roster
\item"$(i)$"
       for some $l\in\{0,1,2\}$, the number of matrices among $A_1,\dots,A_m$
       belonging to $C_3^{(*,l)}$ is odd;
\item"$(ii)$"
       $i_1,\dots,i_m\in\{2,6\}$ and the number of $2$'s in the sequence $(i_1,\dots,i_m)$ is odd.
\endroster
\endproclaim


\head\sectNtwo. Products of conjugacy classes in $GU(2,q^2)$ and $SU(2,q^2)$.
\endhead

Let $G$ (resp. $S$; $PS$) be $GU(2,q^2)$ or $GL(2,q)$ (resp.
$SU(2,q^2)$ or $SL(2,q)$; $PSU(2,q^2)$ or $PSL(2,q)$).
We follow the sign convention from \S\sectConvent.

\subhead\sectNTwoG. Class products in $GU(2,q^2)$ and $GL(2,q)$
\endsubhead
%
We use the notation from [\refEnnola] for conjugacy classes in $G$.
The classes (and the
respective Jordan normal forms) are:
$$
 C_1^{(k)}:\;\left(\smallmatrix\omega^k&0\\0&\omega^k\endsmallmatrix\right),\quad
 C_2^{(k)}:\;\left(\smallmatrix\omega^k&0\\1&\omega^k\endsmallmatrix\right),\quad
 C_3^{(k,l)}:\;\left(\smallmatrix\omega^k&0\\0&\omega^l\endsmallmatrix\right),\quad
 C_4^{(k)}:\left(\smallmatrix\rho^k&0\\0&\rho^{\mp qk}\endsmallmatrix\right).
$$
In the last two cases we have $C_3^{(k,l)}=C_3^{(l,k)}$, $C_4^{(k)}=C_4^{(\mp qk)}$ and
we claim that the matrix is non-scalar, i.~e.,
that $k\ne l$ and $k\not\equiv\mp qk\mod q^2-1$ respectively.
There are four families of irreducible characters:
$\chi_1^{(t)}$,
$\chi_q^{(t)}$ ($0\le t\le q$),
$\chi_{q\mp1}^{(t,u)}$ ($1\le t<u\le q\pm1$),
$\chi_{q\pm1}^{(t)}$ ($1\le t\le q^2$, $t\not\equiv0\mod q\mp1$, 
$\chi_{q\pm1}^{(t)}=\chi_{t+1}^{(\mp qt)}$); see details in [\refEnnola].
We denote the union of all $C_i^{(\dots)}$ by $C_i$. We define
$\delta_{a_1,\dots,a_m}$ in the same way as in \S\sectNotation.

\proclaim{ Theorem \thGUtwo }
Let $A_1,\dots,A_m\in G\setminus C_1$, $m\ge 3$, be matrices which satisfy {\rm(\eqDet)}.
Let $A_\nu\in C_{i_\nu}$, $\nu=1,\dots,m$. 
Let 
$$
     i_0=\cases 3,&G=GU,\\ 4,&G=GL,\endcases
	\qquad\text{and}\qquad
     C=\cases C_3^{(0,2)}\cup C_3^{(1,3)}, &G=GU(2,q^2),\\
             C_4^{(2)}, & G=GL(2,q).\endcases
$$
Then $I\not\in A_1^G\dots A_m^G$
if and only if one of the following conditions holds up to permutation of $A_1,\dots,A_m$:
\roster
\item"$(i)$"
             $(i_1,\dots,i_m)=(i_0,i_0,2)$ and $\delta_{111}+\delta_{121}=1$;
\item"$(ii)$"
             $q=3$, $A_1,\dots,A_{m-1}\in C$, 
             and $A_m\in C_2$.
\item"$(iii)$"
             $q=2$, and $2$ occurs an odd number of times in $(i_1,\dots,i_m)$.
\endroster
\endproclaim

\demo{ Proof } Case $m=3$. It is enough to compute the structure constants.
They are listed in Table \tabGUtwo. 

\midinsert
\noindent Table \tabGUtwo. Structure constants for $G=GU(2,q^2)$ or $GL(2,q)$, $A_\nu\in C_{i_\nu}$
\smallskip
\vbox{%
\vbox{\offinterlineskip
\hrule
\halign{&\vrule#&\strut\quad#\hfill\cr
height3pt&\omit&\omit&\omit&&\omit&\omit&\omit&\cr
& $\!\!(i_1,i_2,i_3)$ &\omit& $N_G(A_1,A_2,A_3)/|A_1^G|$ &
& $\!\!(i_1,i_2,i_3)$ &\omit& $N_G(A_1,A_2,A_3)/|A_1^G|$ &\cr
height3pt&\omit&\omit&\omit&&\omit&\omit&\omit&\cr
\noalign{\hrule}
height3pt&\omit&\omit&\omit&&\omit&\omit&\omit&\cr
& $(2,2,2)$ &\omit& $q-2\delta_{111}$ && $(4,3,2)$ &\omit& $q\mp1$ &\cr
& $(3,2,2)$ &\omit& $q\pm1$           && $(4,3,3)$ &\omit& $q\mp1$ &\cr
& $(3,3,2)$ &\omit& $(q\pm1)(1\mp(\delta_{111}+\delta_{121}))\;$ &
& $(4,4,2)$ &\omit& $(q\mp1)(1\pm(\delta_{111}+\delta_{121}))\;$ &\cr
& $(3,3,3)$ &\omit& $q\pm1 \mp q\Delta$ && $(4,4,3)$ &\omit& $q\mp1$             &\cr
& $(4,2,2)$ &\omit& $q\mp1$             && $(4,4,4)$ &\omit& $q\mp1 \pm q\Delta$ &\cr
height3pt&\omit&\omit&\omit&&\omit&\omit&\omit&\cr
}\hrule}
\vskip1pt
$\Delta=\delta_{111}+\delta_{112}+\delta_{121}+\delta_{211}$}
\endinsert

Case $m=4$. Suppose that $q\ge 4$. Let $C'=C_4$ if $G=GU$ and $C'=C_3$ if $G=GL$.
Then for any $d\in\UC$ there exists $B\in C'$ such that $\det B=d$.
Hence we can choose $B\in C'$ such that $\det B=\det(A_1A_2)$. Then it follows from
the above computations for $m=3$ that
$B\in A_1^G A_2^G$ and $B^{-1}\in A_3^G A_4^G$.

When $q=3$, the result easily follows from the following fact. If $(A_1,A_2,A_3)$ is a
triple of non-scalar matrices which does not satisfy $(ii)$, then $A_1^G A_2^G A_3^G$
is a coset of $S$ in $G$, maybe, with one scalar matrix missing.
If $q=2$, then $G$ is isomorphic to $S_3\times\UC$.
\qed\enddemo

\subhead\sectNTwoCCS. Conjugacy classes in $SU(2,q^2)\cong SL(2,q)$
\endsubhead
In this section we do not apply the convention of \S\sectConvent.
We use here ``$SU$-language'' but, using Table \tabSUtwoSLtwo,
everything can be easily translated to ``$SL$-language''.
So, we set $S=SU(2,q^2)$ and $G=GU(3,q^2)$ and the notation $C_i^{(\dots)}$ is used
for conjugacy classes of $G$ and $S$ (except the second column of Table \tabSUtwoSLtwo).

It is known that $S$ is isomorphic to $SL(2,q)$.
In fact, these groups are conjugated in $GL(2,q^2)$ (but not in $SL(2,q^2)$!).
Indeed, let $z\in\GF{q^2}$ be such that $\bar z=-z$. Then the Hermitian form
$\left(\smallmatrix 0&z\\ \bar z& 0\endsmallmatrix\right)$ is preserved by any element
of $SL(2,q)$. We fix an isomorphism $\Phi:SU(2,q^2)\to SL(2,q)$.

If $q$ is even, then $G=S\times\UC$, so the class product problem for
$S$ is reduced to that for $G$ (see \S\sectNtwoP\ for more details).
So, we suppose that $q=p^m=2r-1$.
We set also $r'=r-1$ (so, $q=2r'+1$). In this case we can choose $z=\rho^r$.

The conjugacy classes of $S$ are as follows.
Each of $C_2^{(k)}$, $k=0,r$, splits into two classes
$C_2^{(k,l)}$, $l=0,1$ so that
$\Phi\big(C_2^{(k,l)}\big)$ is the conjugacy class in $SL(2,q)$ of
$(-1)^{k/r} \left(\smallmatrix1&0\\ \sigma^l &1\endsmallmatrix\right)$
where $\sigma=\rho^{q+1}$ is a generator of $\GF{q}^*$.
This notation of conjugacy classes in $S$ depends on the choice of $\Phi$.

Other conjugacy classes of $G$ contained in $S$ are conjugacy
classes of $S$. The list of all conjugacy classes of the both groups and the
correspondence between them under the isomorphism $\Phi$ is given in Table
\tabSUtwoSLtwo.

\midinsert
\noindent Table \tabSUtwoSLtwo. Correspondence of classes in $SU(2,q^2)$ and $SL(2,q)$, $q=3r-1=3r'+1$
\smallskip
\vbox{\offinterlineskip
\def\TA{Range of the parameters}
\def\TB{$k=0,1$; $\;\;l=0,1$}
\def\TC{$(q+1)/\gcd(q+1,k)$}
\def\TD{$(q-1)/\gcd(q-1,k)\;\;\;$}
\hrule
\halign{&\vrule#&\strut\quad#\hfill\cr
height3pt&\omit    &\omit& \omit            &\omit& \omit            &\omit& \omit    &\cr
& Class in $SU$    &\omit& Class in $SL$    &\omit& \TA              &\omit& Order    &\cr
height3pt&\omit    &\omit& \omit            &\omit& \omit            &\omit& \omit    &\cr
\noalign{\hrule}
height3pt&\omit    &\omit& \omit            &\omit& \omit            &\omit& \omit    &\cr
& $C_1^{(rk)}$     &\omit& $C_1^{(r'k)}$    &\omit& $k=0,1$          &\omit& $k+1$    &\cr
& $C_2^{(rk,l)}$   &\omit& $C_2^{(r'k,l)}$  &\omit& \TB              &\omit& $(k+1)p$ &\cr
& $C_3^{(k,-k)}$   &\omit& $C_4^{((q-1)k)}$ &\omit& $k=1,\dots,r-1$  &\omit& \TC      &\cr
& $C_4^{((q+1)k)}$ &\omit& $C_3^{(k,-k)}$   &\omit& $k=1,\dots,r'-1$ &\omit& \TD      &\cr
height3pt&\omit    &\omit& \omit            &\omit& \omit            &\omit& \omit    &\cr
}\hrule}
\endinsert

The class product problem for pairs of matrices
(to determine the class of the inverse matrix) has an evident solution for
$C_1,C_3,C_4$. The answer for $C_2$ is:

\proclaim{ Proposition \propSUtwoInv }
Let $A\in SU(2,q^2)$, $q=2r-1$. Let $A\in C_2^{(k,l)}$, $k,rl\in\{0,r\}$.
Then
$A^{-1}\in C_2^{(k,l)}$ when $r$ is odd and
$A^{-1}\in C_2^{(k,1-l)}$ when $r$ is even.
\endproclaim

\demo{ Proof }
This follows from the fact that two matrices
$\left(\smallmatrix 1&0\\a&1\endsmallmatrix\right)$ and
$\left(\smallmatrix 1&0\\b&1\endsmallmatrix\right)$, $ab\ne 0$ are
conjugated in $SL(2,K)$ if and only if $ab$ is a square in $K$.
\qed\enddemo

\remark{ Remark \remSUtwo } (cp. Remark in \S\sectCCS).
Let $C$ be the conjugacy class of $\left(\smallmatrix1&0\\1&1\endsmallmatrix\right)$ in $GL(2,q^2)$.
Then $C\cap SL(2,q^2)$ splits into two classes, let us denote them by $C^{(0)}$ and $C^{(1)}$ 
However, the splitting of $C_2^{(0)}$ in $SU$ does not follow the splitting of $C$.
We have $C^{(1)}\cap SU=\varnothing$ and $C^{(0)}\cap SU=C_2$.
This is why there is no any canonical form of these classes in $SL(2,q^2)$.
\endremark

\bigskip

\subhead\sectNTwoS. Class products in $SU(2,q^2)\cong SL(2,q)$
\endsubhead

\proclaim{ Theorem \thSUtwo }
Let $G=GU(2,q^2)$, $S=SU(2,q^2)$, $q=2r-1$. 
Let $A_1,\dots,A_m\in S\setminus C_1$, $m\ge 3$, be such that
$I\in A_1^G\dots A_m^G$. Then $I\not\in A_1^S\dots A_m^S$ if and only if $m=3$ and
one of the following conditions holds up to change of the order of $A_1,\dots,A_m$:
\roster
\item"$(i)$"
        $m=3$, $A_\nu\in C_2^{(rk_\nu,l_\nu)}$ $(\nu=1,2,3)$,
        $l_1\ne l_2$, and $\delta_{111}=0$ {\rm(}i.~e., $k_1+k_2+k_3$ is odd{\rm);}
\item"$(ii)$"
        $m=3$, $A_\nu\in C_2^{(rk_\nu,l_\nu)}$ $(\nu=1,2)$,
        $A_3\in C_3^{(k_3,-k_3)}$, and $r(k_1+k_2+1)+k_3+l_1+l_2$ is odd
        {\rm(see Table \tabSUtwo);}
\item"$(iii)$"
        $m=3$, $A_\nu\in C_2^{(k_\nu,l_\nu)}$ $(\nu=1,2)$,
        $A_3\in C_4^{((q+1)k_3)}$, and $(r-1)(k_1+k_2+1)+k_3+l_1+l_2$ is odd
        {\rm(see Table \tabSUtwo);}
\item"$(iv)$"
        $q=3$ and $\varphi(A_1)+\dots+\varphi(A_m)\not\equiv0\mod3$ where
        $\varphi(A)=l+1$ if $A\in C_2^{(2k,l)}$ and
        $\varphi(A)=0$ if $A\not\in C_2$ {\rm(see Remark \remSUtwoF\rm);}
\item"$(v)$"
        $q=5$, $m=3$, $A_\nu\in C_2^{(3k_\nu,l_\nu)}$ $(\nu=1,2,3)$,
        $l_1=l_2=l_3$, and
        $\delta_{111}=1$ {\rm(}i.~e., $k_1+k_2+k_3$ is even{\rm);}
\item"$(vi)$"
        $q=5$, $m=4$, $A_\nu\in C_2^{(3k_\nu,l_\nu)}$ $(\nu=1,2,3,4)$,
        $l_1=\dots=l_4$, and $\delta_{1111}=0$ {\rm(}i.~e., $k_1+\dots+k_4$ is odd{\rm);}
\item"$(vii)$"
        $q=5$, $m=4$, $A_\nu\in C_2^{(3k_\nu,l_\nu)}$ $(\nu=1,2,3)$, $A_4\in C_3^{(k_4,-k_4)}$,
        $l_1=l_2=l_3$, and $k_1+\dots+k_4$ is odd.
\endroster
\endproclaim

\remark{ Remark \remSUtwoF } The mapping
$a\mapsto\left(\smallmatrix1& 0\\1&1\endsmallmatrix\right)$,
$i\mapsto\left(\smallmatrix0&-1\\1&0\endsmallmatrix\right)$
defines an isomorphism $F\cong SL(2,3)$ where $F$ is the group
discussed in \S\sectSmallqG.
The class products in $F$ are described in Proposition \propSmallqOne.
The correspondence of classes
between $F$ and $SU(2,3^2)$ is:
$\{(-1)^k\}\to C_1^{(2k)}$,
$i^F\to C_3^{(1,-1)}$,
$(-1)^k a^F\to C_2^{(2k,0)}$,
$(-1)^k b^F\to C_2^{(2k,1)}$.
\endremark


\midinsert
\noindent Table \tabSUtwo. ${N_{SU(2,q^2)}(A_1,A_2,A_3)\over qr(q-1)}$ for $q=2r-1$,
$A_1\in C_2^{(rk_1,l_1)}$, $A_2\in C_2^{(rk_2,l_2)}$, $A_3\in C_3\cup C_4$
\smallskip
\vbox{\offinterlineskip
\hrule
\halign{&\vrule#&\strut\hfill#\hfill\cr
height3pt
&\omit&
&\multispan{15}&&\multispan{15}&\cr
&\;$r$\;&
&\multispan{15}\hfill $r$ even \hfill&&\multispan{15}\hfill $r$ odd \hfill&\cr
height3pt
&\omit&
&\multispan{15}&&\multispan{15}&\cr
\noalign{\hrule}
height3pt
&\omit&
&\multispan 7&&\multispan 7&&\multispan 7&&\multispan 7&\cr
&\;$l_1,l_2$\;&
&\multispan7\hfill$l_1=l_2$\hfill&&\multispan7\hfill$l_1\ne l_2$\hfill&
&\multispan7\hfill$l_1=l_2$\hfill&&\multispan7\hfill$l_1\ne l_2$\hfill&\cr
height3pt
&\omit&
&\multispan 7&&\multispan 7&&\multispan 7&&\multispan 7&\cr
\noalign{\hrule}
height3pt
&\omit&
&\multispan 3&&\multispan 3&&\multispan 3&&\multispan 3&
&\multispan 3&&\multispan 3&&\multispan 3&&\multispan 3&\cr
&\;\;\;\;$k_1+k_2\mod 2$\;\;\;\;&
&\multispan 3\hfill 0\hfill&&\multispan 3\hfill$1$\hfill&
&\multispan 3\hfill 0\hfill&&\multispan 3\hfill$1$\hfill&
&\multispan 3\hfill 0\hfill&&\multispan 3\hfill$1$\hfill&
&\multispan 3\hfill 0\hfill&&\multispan 3\hfill$1$\hfill&\cr
height3pt
&\omit&
&\multispan 3&&\multispan 3&&\multispan 3&&\multispan 3&
&\multispan 3&&\multispan 3&&\multispan 3&&\multispan 3&\cr
\noalign{\hrule}
height3pt
&\omit&
&\omit&&\omit&&\omit&&\omit&&\omit&&\omit&&\omit&&\omit&
&\omit&&\omit&&\omit&&\omit&&\omit&&\omit&&\omit&&\omit&\cr
&\;$k_3\mod 2$\;&
&\;\;0\;\;&&\;\;1\;\;&&\;\;0\;\;&&\;\;1\;\;&&\;\;0\;\;&&\;\;1\;\;&&\;\;0\;\;&&\;\;1\;\;&
&\;\;0\;\;&&\;\;1\;\;&&\;\;0\;\;&&\;\;1\;\;&&\;\;0\;\;&&\;\;1\;\;&&\;\;0\;\;&&\;\;1\;\;&
\cr
height3pt
&\omit&
&\omit&&\omit&&\omit&&\omit&&\omit&&\omit&&\omit&&\omit&
&\omit&&\omit&&\omit&&\omit&&\omit&&\omit&&\omit&&\omit&\cr
\noalign{\hrule}
height3pt
&\multispan{33}&\cr
\noalign{\hrule}
height3pt
&\omit&
&\omit&&\omit&&\omit&&\omit&&\omit&&\omit&&\omit&&\omit&
&\omit&&\omit&&\omit&&\omit&&\omit&&\omit&&\omit&&\omit&\cr
&\;$A_3\in C_3^{(k_3,-k_3)}$\;&
&\;1\;&&\;0\;&&\;1\;&&\;0\;&&\;0\;&&\;1\;&&\;0\;&&\;1\;&
&\;0\;&&\;1\;&&\;1\;&&\;0\;&&\;1\;&&\;0\;&&\;0\;&&\;1\;&
\cr
height3pt
&\omit&&\omit&&\omit&&\omit&&\omit&&\omit&&\omit&&\omit&&\omit&
&\omit&&\omit&&\omit&&\omit&&\omit&&\omit&&\omit&&\omit&\cr
\noalign{\hrule}
height3pt
&\omit&
&\omit&&\omit&&\omit&&\omit&&\omit&&\omit&&\omit&&\omit&
&\omit&&\omit&&\omit&&\omit&&\omit&&\omit&&\omit&&\omit&\cr
&\;$A_4\in C_4^{((q+1)k_3)}$\;&
&\;0\;&&\;1\;&&\;1\;&&\;0\;&&\;1\;&&\;0\;&&\;0\;&&\;1\;&
&\;1\;&&\;0\;&&\;1\;&&\;0\;&&\;0\;&&\;1\;&&\;0\;&&\;1\;&
\cr
height3pt
&\omit&&\omit&&\omit&&\omit&&\omit&&\omit&&\omit&&\omit&&\omit&
&\omit&&\omit&&\omit&&\omit&&\omit&&\omit&&\omit&&\omit&\cr
}\hrule}
\endinsert

\demo{ Proof }
Case $m=3$.
It is enough to consider only the triples $(A_1,A_2,A_3)$ containing
at least two matrices from $C_2$ (otherwise $N_S(A_1,A_2,A_3)=N_G(A_1,A_2,A_3)$).
We compute $N_S(A_1,A_2,A_3)$ for all such triples. If $A_\nu\in C_2^{(k_\nu,l_\nu)}$,
$\nu=1,2,3$, then we have
$$
   N_S(A_1,A_2,A_3)=\cases r(r-1)(2q-(3r-3e_r+1)\delta_{111}),& l_1=l_2=l_3,\\
                  r(r-1)(r-e_r-1)\delta_{111},& l_1=l_2\ne l_3,
   \endcases
$$
where $e_r={1+(-1)^r\over 2}$.

If $A_1\in C_2^{(rk_1,l_1)}$, $A_2\in C_2^{(rk_2,l_2)}$, and $A_3\in C_3\cup C_4$, we have
$N_S(A_1,A_2,A_3)=qr(q-1)\delta^*$ where the values of $\delta^*$ are given in Table \tabSUtwo.

\medskip
Case $m\ge 4$. The result for $m>4$ follows from the result for $m=4$. So we assume that $m=4$.
If $q=3$, then $S$ is isomorphic to the group $F$ discussed in \S\sectSmallq\ and the result
follows from Proposition \propSmallqOne.
If $q=5$, then $S$ then it is enough to compute explicitly the structure constants for
all triples and quadruples. So, we assume that $q\ge 7$.

If one of $A_\nu$ does not belong to $C_2$, then we can choose $B\in C_4$ such that
$B\in A_1^S A_2^S$ and $B^{-1}\in A_3^SA_4^S$.

If $A_\nu\in C_2^{(rk_\nu,l_\nu)}$, $k=1,\dots,4$, then without loss of generality we
may assume that $l_3=l_4=l$. Let $B\in C_2^{(rk,l)}$ where $k+k_1+k_2$ is even.
Then $B\in A_1^S A_2^S$ and $B^{-1}\in A_3^S A_4^S$.
\qed\enddemo


\bigskip

\subhead\sectNTwoP. Class products in $PSU(2,q^2)\cong PSL(2,q)$
\endsubhead

Let $PS=PSU(2,q^2)\cong PSL(2,q)$, $q\ge4$.
Like in Corollary \corPSU, we denote the projection of a class $C_i^{(\dots)}$ by
$\tilde C_i^{(\dots)}$. Products of conjugacy classes in $PS$ are partially computed in
[\refAH; Ch.~4, Th.~4.2]. For reader conenience we give the correspondence of
notation in Tables \tabPSUtwo.1 -- \tabPSUtwo.2.

\midinsert
\noindent Table \tabPSUtwo.1. Conjugacy classes in $SU(2,q^2)=PSU(2,q^2)\cong SL(2,q)=PSL(2,q)$
for even $q$.
\smallskip
\vbox{\offinterlineskip
\def\TA{Parameters}
\def\TB{$C_4^{((q-1)k)}$}
\def\TC{$k=1,\dots,{q\over2}$}   \def\TD{$(q+1)/\gcd(q+1,k)$}
\def\TE{$k=1,\dots,{q-2\over2}$} \def\TF{$(q-1)/\gcd(q-1,k)\quad$}
\hrule
\halign{&\vrule#&\strut\quad#\hfill\cr
height3pt&\omit &\omit& \omit            &\omit& \omit          &\omit& \omit &\omit& \omit &\cr
& In [\refAH]   &\omit& Class in $SU$    &\omit& Class in $SL$  &\omit& \TA   &\omit& Order &\cr
height3pt&\omit &\omit& \omit            &\omit& \omit          &\omit& \omit &\omit& \omit &\cr
\noalign{\hrule}
height3pt&\omit &\omit& \omit            &\omit& \omit          &\omit& \omit &\omit& \omit &\cr
& $C_1$         &\omit& $C_1^{(0)}$      &\omit& $C_1^{(0)}$    &\omit& \omit &\omit& $1$   &\cr
& $C_2$         &\omit& $C_2^{(0)}$      &\omit& $C_2^{(0)}$    &\omit& \omit &\omit& $2$   &\cr 
& $R_k$         &\omit& $C_3^{(k,-k)}$   &\omit& \TB            &\omit& \TC   &\omit& \TD   &\cr
& $K_k$         &\omit& $C_4^{((q+1)k)}$ &\omit& $C_3^{(k,-k)}$ &\omit& \TE   &\omit& \TF   &\cr
height3pt&\omit &\omit& \omit            &\omit& \omit          &\omit& \omit &\omit& \omit &\cr
}\hrule}
\endinsert

\midinsert
\noindent Table \tabPSUtwo.2. Conjugacy classes in $PSU(2,q^2)\cong PSL(2,q)$
for $q=p^m=2r-1$ (for a prime $p$), $r'=r-1$.
\smallskip
\vbox{\offinterlineskip
\def\TA{Parameters}
\def\TB{$\tC_4^{((q-1)k)}$}
\def\TC{$k=1,\dots,[{r\over2}]$}   \def\TD{$r/\gcd(r,k)$}
\def\TE{$k=1,\dots,[{r'\over2}]$} \def\TF{$r'/\gcd(r',k)\quad$}
\def\tC{\tilde C}
\hrule
\halign{&\vrule#&\strut\quad#\hfill\cr
height3pt&\omit &\omit& \omit             &\omit& \omit            &\omit& \omit &\omit& \omit &\cr
& In [\refAH]   &\omit& Class in $PSU$    &\omit& Class in $PSL$   &\omit& \TA   &\omit& Order &\cr
height3pt&\omit &\omit& \omit             &\omit& \omit            &\omit& \omit &\omit& \omit &\cr
\noalign{\hrule}
height3pt&\omit &\omit& \omit             &\omit& \omit            &\omit& \omit &\omit& \omit &\cr
& $C_1$         &\omit& $\tC_1^{(0)}$     &\omit& $\tC_1^{(0)}$    &\omit& \omit &\omit& $1$   &\cr
& $C_2$         &\omit& $\tC_2^{(0,0)}$   &\omit& $\tC_2^{(0,0)}$  &\omit& \omit &\omit& $p$   &\cr 
& $C_3$         &\omit& $\tC_2^{(0,1)}$   &\omit& $\tC_2^{(0,1)}$  &\omit& \omit &\omit& $p$   &\cr 
& $R_k$         &\omit& $\tC_3^{(k,-k)}$  &\omit& \TB              &\omit& \TC   &\omit& \TD   &\cr
& $K_k$         &\omit& $\tC_4^{((q+1)k)}$&\omit& $\tC_3^{(k,-k)}$ &\omit& \TE   &\omit& \TF   &\cr
height3pt&\omit &\omit& \omit             &\omit& \omit            &\omit& \omit &\omit& \omit &\cr
}\hrule}
\endinsert

As in the previous section, we use here the ``$SU$-notation'' for conjugacy classes in $PS$
(the second column in Tables \tabPSUtwo.1 -- \tabPSUtwo.2).

\proclaim{ Corollary \corPSUtwo } Let $m\ge3$, $q\ge4$, and $c_1,\dots,c_m$ are non-identity
conjugacy classes in $PS$.
Then $I\not\in c_1\dots c_m$ if and only if $m=3$ and one of the following cases occurs up
to permutation:
\roster
\item "$(i)$"
$q$ is even and $(c_1,c_2,c_3)=\big(\,C_3^{(k,-k)},\,C_3^{(k,-k)},\,C_2^{(0)}\,\big)$, 
$k=1,\dots,q/2$;
\item "$(ii)$"
$q=2r-1\equiv 1\mod4$ {\rm(}so, $r$ is odd\/{\rm)} and
$(c_1,c_2,c_3)$ is

\smallskip
$\big(\tilde C_2^{(0,l_1)},\,\tilde C_2^{(0,l_2)},\,\tilde C_4^{((q+1)k)}\big)$,
$k+l_1+l_2$ is odd, $k=1,\dots,{r-1\over2}$;

\smallskip
\item "$(iii)$"
$q=2r-1\equiv 3\mod4$ {\rm(}so, $r$ is even\/{\rm)} and $(c_1,c_2,c_3)$ is one of:

$\big(\tilde C_2^{(0,l_1)},\,\tilde C_2^{(0,l_2)},\,\tilde C_3^{(k,-k)}\big)$,
$k+l_1+l_2$ is odd, $k=1,\dots,{r\over2}$;

$\big(\tilde C_2^{(0,l)},\,\tilde C_3^{({r/2},-{r/2})},\,
\tilde C_3^{({r/2},-{r/2})}\big)$, $l=0,1$.
\endroster
\endproclaim

In particular, we see that $\cn(PS)=3$, $\ecn(PS)=4$ (see \S\sectCN).
This fact was already proved in [\refAH; Ch.~4].


\Refs
\def\r{\ref}

\r\no\refAH
\by    Z.~Arad, M.~Herzog (eds.)
\book  Products of conjugacy classes in groups \bookinfo Lecture Notes in Math. 1112
\publ  Springer-Verlag  \publaddr Berlin, Heidelberg, N.Y., Tokyo \yr 1985
\endref

\r\no\refAW
\by     S.~Agnihotri, C.~Woodward
\paper  Eigenvalues of products of unitary matrices and quantum Schubert
        calculus \jour Math. Research Letters \vol 5 \yr 1998 \pages 817--836
\endref

\r\no\refBelkale
\by     P.~Belkale
\paper  Local systems on $\Bbb P^1-S$ for $S$ a finite set
        \jour Compos. Math. \vol 129 \yr 2001 \pages 67--86
\endref

\r\no\refDickson
\by     L.~E.~Dickson
\book   Linear groups with an exposition of the Galois field theory
\publ   Teubner \publaddr Leipzig \yr 1901
\endref

\r\no\refEnnolaConj
\by     V.~Ennola
\paper  On the conjugacy classes of the finite unitary groups
\jour   Ann. Acad. Sci. Fennicae, Ser. A I.  \vol 313\yr 1962 \pages 3--13
\endref

\r\no\refEnnola
\by     V.~Ennola
\paper  On the characters of the finite unitary groups
\jour   Ann. Acad. Sci. Fennicae, Ser. A I.  \vol 323\yr 1963 \pages 3--35
\endref

\r\no\refGAP
\book   GAP software; the file {\tt ctgeneri.tbl}
\endref

\r\no\refGeck
\by     M.~Geck \paper Diploma thesis
\endref

\r\no\refGordeev
\by     N.~L.~Gordeev
\paper  Products of conjugacy classes in perfect linear groups. Extended covering number
\jour   Zapiski Nauchn. Semin. POMI \vol 321 \yr 2005 \pages 67--89 \lang Russian
\transl English transl.
\jour   J. Math. Sci. \vol 136 \yr 2006 \pages 3867--3879
\endref

\r\no\refKarni
\by     S.~Karni
\paper  Covering number of groups of small order and sporadic groups
\inbook Ch.~3 in [\refAH] \pages 52--196
\endref

\r\no\refLev
\by     A.~Lev
\paper  The covering number of the group $\operatorname{PSL}_n(F)$
\jour   J. Algebra \vol 182 \yr 1996 \pages 60--84
\endref

\r\no\refLSh
\by     M.~W.~Liebeck, A.~Shalev
\paper  Diameter of finite simple groups: sharp bounds and applications
\jour   Ann. of Math. \vol 154 \yr 2001 \pages 383--406
\endref

\r\no\refMM
\by     G.~Malle, B.~H.~Matzat
\book   Inverse Galois Theory \publ Springer-Verlag \publaddr Berlin, N.Y. \yr 1999
\endref

\r\no\refOrJKTR
\by     S.~Yu.~Orevkov
\paper  Quasipositivity test via unitary representations of braid groups and its
        applications to real algebraic curves
\jour   J. of Knot Theory and Ramifications  \vol 10 \yr 2001 \pages 1005--1023
\endref

\r\no\refSF
\by     W.A.~Simpson, J.S.~Frame
\paper  The character tables for
        $\operatorname{SL}(3,q)$,
        $\operatorname{SU}(3,q^2)$,
        $\operatorname{PSL}(3,q)$,
        $\operatorname{PSU}(3,q^2)$
\jour Canad. J. Math. \vol 25\yr 1973 \pages 486--494
\endref

\r\no\refWallG
\by     G.~E.~Wall
\paper  On the conjugacy classes in the unitary, symplectic and orthogonal groups
\jour   J. Austral. Math. Soc. \vol 3 \yr 1963 \pages 1--62
\endref

\r\no\refWallSL
\by     G.~E.~Wall
\paper  Conjugacy classes in projective and special linear groups
\jour   Bull. Austral. Math. Soc. \vol 22 \yr 1980 \pages 339-364
\endref

\endRefs

\enddocument